%%%%%%%%%%%%%%%%%%%%%%%%%%%%%%%%
%
% These are the macroes
%
%%%%%%%%%%%%%%%%%%%%%%%%%%%%%%%%

\magnification\magstephalf

\voffset0truecm
\hoffset=0truecm
\vsize=23truecm
\hsize=15.8truecm
\topskip=1truecm

\binoppenalty=10000
\relpenalty=10000

%%%% Famiglie di Fonts
\font\tenbb=msbm10		\font\sevenbb=msbm7		\font\fivebb=msbm5
\font\tensc=cmcsc10		\font\sevensc=cmcsc7 	\font\fivesc=cmcsc5
\font\tensf=cmss10		\font\sevensf=cmss7		\font\fivesf=cmss5
\font\tenfr=eufm10		\font\sevenfr=eufm7		\font\fivefr=eufm5

%\font\eurb=eurb10

\newfam\bbfam	\newfam\scfam	\newfam\frfam	\newfam\sffam

\textfont\bbfam=\tenbb
\scriptfont\bbfam=\sevenbb
\scriptscriptfont\bbfam=\fivebb

\textfont\scfam=\tensc
\scriptfont\scfam=\sevensc
\scriptscriptfont\scfam=\fivesc

\textfont\frfam=\tenfr
\scriptfont\frfam=\sevenfr
\scriptscriptfont\frfam=\fivefr

\textfont\sffam=\tensf
\scriptfont\sffam=\sevensf
\scriptscriptfont\sffam=\fivesf

\def\bb{\fam\bbfam \tenbb} % Blackboard bold
\def\sc{\fam\scfam \tensc} % Maiuscoletto
\def\fr{\fam\frfam \tenfr} % Gotico
 % Sans serif

%%%% Fonts

\font\sezfont=cmbx10 scaled \magstep1
\font\subsectfont=cmbx10 scaled \magstephalf
\font\titfont=cmbx10 scaled \magstep2
\font\autfont=cmcsc10
\font\intfont=cmss10 

%%%% Abbreviazioni
\let\no=\noindent
\let\bi=\bigskip
\let\me=\medskip
\let\sm=\smallskip
\let\ce=\centerline

\let\io=\infty

%%%% Numerazione

\newcount\sectno\sectno=0
\newcount\subsectno\subsectno=0
\newcount\thmno\thmno=0
\newcount\tagno\tagno=0
\newcount\notitolo\notitolo=0
\newcount\defno\defno=0

% sezione
\def\sect#1\par{
	\global\advance\sectno by 1 \global\subsectno=0\global\defno=0\global\thmno=0
	\vbox{\vskip.75truecm\advance\hsize by 1mm
	\hbox{\centerline{\sezfont \the\sectno.~~#1}}
	\vskip.25truecm}\nobreak}

% sottosezione
\def\subsect#1\par{
	\global\advance\subsectno by 1
	\vbox{\vskip.75truecm\advance\hsize by 1mm
	\line{\subsectfont \the\sectno.\the\subsectno~~#1\hfill}
	\vskip.25truecm}\nobreak}
	
% definizioni
\def\defin#1{\global\advance\defno by 1
	\global\expandafter\edef\csname+#1\endcsname%
    {\number\sectno.\number\defno}
    \no{\bf Definition~\the\sectno.\the\defno.}}

% teoremi, lemmi, proposizioni e osservazioni
\def\thm#1#2{
	\global\advance\thmno by 1
	\global\expandafter\edef\csname+#1\endcsname%
	{\number\sectno.\number\thmno}
	\no{\bf #2~\the\sectno.\the\thmno.}}

% equazioni
\def\Tag#1{\global\advance\tagno by 1 {(\the\tagno)}
    \global\expandafter\edef\csname+#1\endcsname%
    		{(\number\tagno)}}
\def\tag#1{\leqno\Tag{#1}}

\def\rf#1{\csname+#1\endcsname\relax}

% fine dimostrazione
\def\proof{\no{\sl Proof.}\enskip}
\def\qedn{\thinspace\null\nobreak\hfill\hbox{\vbox{\kern-.2pt\hrule height.2pt
        depth.2pt\kern-.2pt\kern-.2pt \hbox to2.5mm{\kern-.2pt\vrule
        width.4pt \kern-.2pt\raise2.5mm\vbox to.2pt{}\lower0pt\vtop
        to.2pt{}\hfil\kern-.2pt \vrule
        width.4pt\kern-.2pt}\kern-.2pt\kern-.2pt\hrule height.2pt
        depth.2pt \kern-.2pt}}\par\medbreak}
    \def\qed{\hfill\qedn}
    
%%% Numerazione di pagina / intestazioni
\newif\ifpage\pagefalse
\newif\ifcen\centrue

\headline={
\ifcen\hfil\else
\ifodd\pageno
\global\hoffset=0.5truecm
\else
\global\hoffset=-0.4truecm
\fi\hfil
\fi}

\footline={
	\ifpage
		\hfill\rm\folio\hfill
	\else
		\global\pagetrue\hfill
\fi}

\lccode`\'=`\'

%%% Riferimenti bibliografici
\def\bib#1{\me\item{[#1]\enskip}}

%%% Macrodefinizioni matematiche
\def\ca#1{{\cal #1}}
\def\C{{\bb C}} \def\d{{\rm d}}
\def\R{{\bb R}} \def\Q{{\bb Q}}
\def\Z{{\bb Z}} \def\N{{\bb N}}
\def\T{{\bb T}}

\def\X{{\fr X}} \def\S{{\bb S}}
\let\la=\langle	\let\ra=\rangle
\let\eps=\varepsilon \let\phe=\varphi
\let\de=\partial

\mathchardef\void="083F
 % PARTE REALE
\def\Im{\mathop{\rm Im}\nolimits} % PARTE IMMAGINARIA

\def\Aut{\mathop{\rm Aut}\nolimits}

\def\diag{\mathop{\rm Diag}\nolimits}

\def\id{\mathop{\rm Id}\nolimits}

\def\s{{\rm s}}

\def\res{{\rm Res}}

\def\invlim{\mathop{\vtop{\offinterlineskip
\hbox{\rm lim}\kern1pt\hbox{\kern-1.5pt$\longleftarrow$}\kern-3pt}
}\limits}
 % "PER OGNI" messi al posto giusto

%%%%%%%%%%%%%%%%%%%%%%%%%%%%%%%%%%%%%%%%%%%%%%%%%%%%%%%%%%%%%%%%%%%%%%%%%%%%%%%%%%%%%%%%%%%%%%%%%%%%%%%%%%%%%%%%%%%%%%%%%%%%%%%%%%%%%%%%%%%%%%%%%%%%%%%%%%%%%%%%%%%%%%%%%%%%%%%%%%%%%%%%%%%%%%%%%%%%%%%%%%%%%%%%%%%%%%%%%%%%%%%%%

%%%%%%%%%%%%%%%%%%%%%%%%%%%%%%%%
%
% Begin of the paper
%
%%%%%%%%%%%%%%%%%%%%%%%%%%%%%%%%

\input graphicx.tex
\input xy
\xyoption{all}

\ce{\titfont Torus actions in the normalization problem} 

\me\ce{\autfont Jasmin Raissy}
\sm\ce{\intfont Dipartimento di Matematica, Universit\`a di Pisa}

\ce{\intfont Largo Bruno Pontecorvo 5, 56127 Pisa}

\sm\ce{\intfont E-mail: {\tt raissy@mail.dm.unipi.it}}
\bi

{\narrower

{\sc Abstract.} Let $f$ be a germ of biholomorphism of $\C^n$, fixing the origin. We show that if the germ commutes with a torus action, then we get information on the germs that can be conjugated to $f$, and furthermore on the existence of a holomorphic linearization or of a holomorphic normalization of $f$. We find out in a complete and computable manner what kind of structure a torus action must have in order to get a Poincar\'e-Dulac holomorphic normalization, studying the possible torsion phenomena. In particular, we link the eigenvalues of $\d f_O$ to the weight matrix of the action. The link and the structure we found are more complicated than what one would expect; a detailed study was needed to completely understand the relations between torus actions, holomorphic Poincar\'e-Dulac normalizations, and torsion phenomena. We end the article giving an example of techniques that can be used to construct torus actions.
 
}

\bi

\footnote{}{{\sl \hskip-20pt\noindent Mathematics Subject Classification (2000).} Primary 37F99; Secondary 32H50, 37G05. \hfill \break {\sl Key words and phrases.} Normalization problem, linearization problem, resonances, torus actions, discrete local holomorphic dynamical systems.}
%{\sf [versione 06/08/2009]}

%%%%%%%%%%%%%%%%%%%%%%%%%%%%%%%%%%%%%%%%%%%%%%%%%%%%%%%%%%%%%%%%%%%%%%%%%%%%%%%%
\sect Introduction
%%%%%%%%%%%%%%%%%%%%%%%%%%%%%%%%%%%%%%%%%%%%%%%%%%%%%%%%%%%%%%%%%%%%%%%%%%%%%%%%

\sm We consider a germ of biholomorphism~$f$ of~$\C^n$ at a fixed point~$p$, which we may place at the origin~$O$. One of the main questions in the study of local holomorphic dynamics (see [A1], [A2], and [Bra] for general surveys on this topic) is when~$f$ is {\it holomorphically linearizable}, i.e., when there exists a local holomorphic change of coordinates such that~$f$ is conjugated to its linear part. The answer to this question depends on the set of eigenvalues of~$\d f_O$, usually called the {\it spectrum} of~$\d f_O$. In fact if we denote by~$\lambda_1, \dots, \lambda_n\in \C^*$ the eigenvalues of~$\d f_O$, then it may happen that there exists a multi-index~$K=(k_1, \dots, k_n)\in \N^n$ with~$|K|=k_1+\cdots+k_n\ge 2$ and such that
$$\lambda^K - \lambda_j:=\lambda_1^{k_1}\cdots\lambda_n^{k_n} - \lambda_j = 0\tag{eqres}$$
for some~$1\le j\le n$; a relation of this kind is called a {\it multiplicative resonance} of~$f$, and~$K$ is called a {\it resonant multi-index}. A {\it resonant monomial} is a monomial~$z^K= z_1^{k_1}\cdots z_n^{k_n}$ in the~$j$-th coordinate such that~$\lambda^K = \lambda_j$. From the formal point of view, we have the following classical result (see [Ar] pp.~192--193 for a proof):

\sm\thm{Te0.1}{Theorem} {\sl Let~$f$ be a germ of holomorphic diffeomorphism of~$\C^n$ fixing the origin~$O$ with no resonances. Then~$f$ is formally conjugated to its differential~$\d f_O$.}

\sm In presence of resonances, even the formal classification is not easy, as the following result of Poincar\'e-Dulac, [Po], [D], shows

\sm\thm{Te0.2}{Theorem} (Poincar\'e-Dulac) {\sl Let~$f$ be a germ of holomorphic diffeomorphism of~$\C^n$ fixing the origin~$O$. Then~$f$ is formally conjugated to a formal power series~$g\in\C[\![z_1, \dots, z_n]\!]^n$ without constant term such that~$\d g_O$ is in Jordan normal form, and~$g$ has only resonant monomials.}

\sm The formal series~$g$ is called a {\it Poincar\'e-Dulac normal form} of~$f$; a proof of Theorem \rf{Te0.2} can be found in [Ar] p.~194. 

The problem with Poincar\'e-Dulac normal forms is that they are not unique. In particular, one may wonder whether it could be possible to have such a normal form including finitely many resonant monomials only. This is indeed the case (see, e.g., Reich [R1]) when $\d f_O$ belongs to the so-called {\it Poincar\'e domain}, that is when $\d f_O$ is invertible and $O$ is either {\it attracting}, i.e., all the eigenvalues of~$\d f_O$ have modulus less than~$1$, or {\it repelling}, i.e., all the eigenvalues of~$\d f_O$ have modulus greater than~$1$ 
(when $\d f_O$ is still invertible but does not belong to the Poincar\'e domain, we shall say that it belongs to the {\it Siegel domain}). 

\sm Even without resonances, the holomorphic linearization is not guaranteed. The easiest positive result is due to Poincar\'e [Po] who, using majorant series, proved the following

\sm\thm{Te0.3}{Theorem} (Poincar\'e, 1893 [Po]) {\sl Let~$f$ be a germ of holomorphic diffeomorphism of~$\C^n$ with an attracting or repelling fixed point. Then~$f$ is holomorphically linearizable if and only if it is formally linearizable. In particular, if there are no resonances then~$f$ is holomorphically linearizable.}

\sm When~$O$ is not attracting or repelling, even without resonances, the formal linearization might diverge. In [Ra2] we found, under certain arithmetic conditions on the eigenvalues and some restrictions on the resonances, a necessary and sufficient condition for holomorphic linearization in presence of resonances, that in fact has as corollaries most of the known linearization results. In [Ra3] we found that, given~$m\ge 2$ germs~$f_1, \dots, f_m$ of biholomorphisms of~$\C^n$, fixing the origin, with~$(\d f_1)_O$ diagonalizable and such that~$f_1$ commutes with~$f_h$ for any~$h=2,\dots, m$, under certain arithmetic conditions on the eigenvalues of~$(\d f_1)_O$ and some restrictions on their resonances,~$f_1, \dots, f_m$ are simultaneously holomorphically linearizable if and only if there exists a particular complex manifold invariant under~$f_1, \dots, f_m$. 

In any cases, there are germs not holomorphically linearizable, for instance when $\d f_O$ is not diagonalizable (see also [PM] for related results):

\sm\thm{Yoccoz}{Theorem} (Yoccoz, 1995 [Y]) {\sl Let $A\in {\rm GL}(n, \C)$ be an invertible matrix such that its eigenvalues have no resonances and such that its Jordan normal form contains a non-trivial block associated to an eigenvalue 
of modulus one. Then there exists a germ of biholomorphism~$f$ of $\C^n$ fixing the origin, with $\d f_O = A$ which is not holomorphically linearizable.}

\sm Then, since every germ of biholomorphism is formally normalizable, studying the {\it holomorphic normalization problem}, i.e., when there exists a local holomorphic change of coordinates such that~$f$ is conjugated to one of its Poincar\'e-Dulac normal forms, could be very useful to understand the dynamics of non-linearizable germs. 

\sm In [Zu], Zung found that to find a Poincar\'e-Dulac holomorphic normalization for a germ of holomorphic vector field is the same as to find (and linearize) a suitable torus action which preserves the vector field.
Following this idea, we found that commuting with a linearizable germ gives us information on the germs conjugated to a given one, and also on the linearization. More precisely we have the following results (for the definition of weight matrix see section 2).

\sm\thm{Toro}{Theorem} {\sl Let $f$ be a germ of biholomorphism of $\C^n$ fixing the origin~$O$. Then $f$ commutes with a holomorphic effective action on $(\C^n,O)$ of a torus of dimension $1\le r\le n$ with weight matrix $\Theta\in M_{n\times r}(\Z)$ if and only there exists a local holomorphic change of coordinates conjugating $f$ to a germ with linear part in Jordan normal form and containing only $\Theta$-resonant monomials.}

\sm\thm{Linearizzazione}{Theorem} {\sl Let~$f$ be a germ of biholomorphism of~$\C^n$ fixing the origin~$O$. Then~$f$ is holomorphically linearizable if and only if it commutes with a holomorphic effective action on~$(\C^n,O)$ of a torus of dimension~$1\le r\le n$ with weight matrix~$\Theta$ having no resonances.}

\sm We want to face and to solve the following problem: to find out in a clear (and possibly computable) manner what kind of structure a torus action must have in order to get a Poincar\'e-Dulac holomorphic normalization from Theorem \rf{Toro}. In particular, to do so we need to link in a clever way the eigenvalues of $\d f_O$ to the weight matrix of the action. Zung dealt with this problem in the case of holomorphic vector fields (see [Zu]), introducing the notion of {\it toric degree} of a vector field. The following definition is a reformulation of Zung's original one, clearer and more suitable to our needs.

\sm\defin{def0.1} The {\it toric degree} of a germ of holomorphic vector field~$X$ of~$(\C^n, O)$ singular at the origin is the minimum~$r\in\N$ such that the semi-simple part $X^{\rm dia}=\sum_{j=1}^n \phe_j z_j \de_j$ of the first jet of $X$ can be written as linear combination with complex coefficients of $r$ diagonal vector fields with integer coefficients, i.e.,  
$$
X^{\s}= \sum_{k=1}^r \alpha_k Z_k,
$$
where~$\alpha_1, \dots, \alpha_r\in\C^*$ and $Z_k = \sum_{j=1}^n\rho^{(k)}_j z_j \de_j$ with~$\rho^{(k)}\in\Z^n$.
The~$r$-tuple~$Z_1,\dots, Z_k$ is called a {\it $r$-tuple of toric vector fields associated to~$X$}, and the numbers~$\alpha_1, \dots, \alpha_n\in\C$ are a~$r$-tuple of {\it toric coefficients} of the toric $r$-tuple.

\sm Then he found that 

\sm\thm{ZungIntro}{Theorem} (Zung, 2002 [Zu]) {\sl Let~$X$ be a germ of holomorphic vector field~$X$ of~$(\C^n, O)$ singular at the origin, of toric degree~$1\le r\le n$. Then~$X$ admits a holomorphic Poincar\'e-Dulac normalization if and only if there exists a holomorphic effective action on~$(\C^n,O)$ of a torus of dimension~$r$ preserving $X$ and such that the columns of the weight matrix of the action are a~$r$-tuple of toric vectors associated to~$X$.}

\me It is a common thinking that once something can be done with germs of vector fields, i.e., for continuous local dynamical systems, then it can be translated analogously for germs of biholomorphisms, i.e., for discrete local dynamical systems. This is not completely true. At the very least there are torsion phenomena to be considered, preventing a straightforward translation from additive resonances (see below for the definition) to multiplicative resonances, and giving rise to new behaviors. One of our aims is exactly to understand up to which point one can push the analogies between continuous and discrete dynamics in the normalization problem. Following \'Ecalle [\'E], we shall use the
following definition of torsion.

\sm\defin{detorsionintro} Let $\lambda\in(\C^*)^n$. The {\it torsion} of~$\lambda$ is the natural integer~$\tau$ such that
$$
{1\over \tau}2\pi i\Z = (2\pi i \Q)\cap \left( (2\pi i \Z)\bigoplus_{1\le j\le n}(\log(\lambda_j)\Z)\right).
$$

To understand what kind of structure a torus action must have in the case of germs of biholomorphisms to get a result equivalent to Theorem \rf{ZungIntro}, we first need a right notion of toric degree for germs of biholomorphisms, and to link it to the torsion we introduced above. The link and the structure we found are more complicated than what one would expect: torsion is not enough to measure the difference between germs of holomorphic vector fields and germs of biholomorphisms. We therefore need a more detailed study.

\sm Notice that given~$\lambda\in(\C^*)^n$, there is a unique~$[\phe]\in(\C/\Z)^n$ such that~$\lambda=e^{2\pi i [\phe]}$, i.e.,~$\lambda_j = e^{2\pi i [\phe_j]}$ for every~$j=1, \dots, n$. The right definition of toric degree for maps is then the following

\sm\defin{De1.1intro} Let~$[\phe]=([\phe_1],\dots, [\phe_n])\in(\C/\Z)^n$, where $[\,\cdot\,]\colon\C^n\to(\C/\Z)^n$ denote the standard projection. The {\it toric degree of~$[\phe]$} is the minimum~$r\in\N$ such that there exist~$\alpha_1, \dots, \alpha_r\in\C$ and~$\theta^{(1)},\dots, \theta^{(r)}\in\Z^n$ such that
$$
[\phe]= \left[\sum_{k=1}^r \alpha_k \theta^{(k)}\right].\tag{Eq1.1}
$$ 
The~$r$-tuple~$\theta^{(1)},\dots, \theta^{(r)}$ is called a {\it $r$-tuple of toric vectors associated to~$[\phe]$}, and the numbers~$\alpha_1, \dots, \alpha_n\in\C$ are {\it toric coefficients} of the toric $r$-tuple.

\sm\defin{DeResIntro} Given~$\theta\in\C^n$ and~$j\in\{1, \dots, n\}$, we say that a multi-index~$Q\in\N^n$, with~$|Q|=\sum_{h=1}^n q_h\ge 2$, gives an {\it additive resonance relation for~$\theta$ relative to the~$j$-th coordinate} if
$$
\langle Q,\theta\rangle =\sum_{h=1}^n q_h\theta_h = \theta_j
$$ 
and we put
$$
{\rm Res}_j^+(\theta)=\{Q\in\N^n\mid |Q|\ge 2,\langle Q,\theta\rangle = \theta_j\}.
$$
Given $[\phe]\in(\C/\Z)^n$, the set
$$
\res_j([\phe])=\{Q\in\N^n\mid |Q|\ge 2,\left[\langle Q,\phe\rangle - \phe_j\right]=[0]\}
$$
of multiplicative resonances of $[\phe]$ is well-defined and we have $\res_j(\lambda)=\res_j([\phe])$, where $\lambda=e^{2\pi i [\phe]}$.

\me We shall find relations between the additive resonances of toric vectors associated to $[\phe]$ and the multiplicative resonances of $[\phe]$. One of the advantages of the approach we found is that we shall be able to easily compute the multiplicative resonances, passing through the additive resonances of $r$-tuples of toric vectors (see Lemma 7.1).

\me Given~$[\phe]\in(\C/\Z)^n$ of toric degree~$1\le r\le n$, even when the $r$-tuple of toric vectors associated to $[\phe]$ is not unique, we can always say whether the toric coefficients are rationally independent with $1$ or not.

\sm\defin{PureTorsionIntro} Let~$[\phe]\in(\C/\Z)^n$ be of toric degree~$1\le r\le n$. We say that $[\phe]$ is in the {\it torsion-free case}, or simply $[\phe]$ {\it is torsion-free}, if its~$r$-tuples of toric vectors have toric coefficients rationally independent with $1$.

\sm As a first application of our methods, we have the following characterization of the vectors~$\lambda\in(\C^*)^n$ without torsion

\sm\thm{PrNoImpureTorsionFreeIntro}{Theorem} {\sl Let $\lambda= e^{2\pi i [\phe]}\in(\C^*)^n$. Then~$[\phe]$ is torsion-free if and only if the torsion of $\lambda$ is $1$.}

\sm In the torsion case, we can always find a more useful toric $r$-tuple.

\sm\defin{De4.2Intro} Let~$[\phe]=([\phe_1],\dots, [\phe_n])\in(\C/\Z)^n$ be of toric degree~$1\le r\le n$ in the torsion case. We say that a~$r$-tuple~$\eta^{(1)},\dots, \eta^{(r)}$ of toric vectors associated to~$[\phe]$ with rationally dependent with $1$ toric coefficients~$\beta_1, \dots, \beta_r$ is {\it reduced} if $\beta_1=1/m$ with~$m\in\N\setminus\{0,1\}$ and~$m, \eta^{(1)}_1, \dots, \eta^{(1)}_n$ coprime. In this case the toric vectors $\eta^{(2)},\ldots,\eta^{(r)}$ are called {\it reduced torsion-free toric vectors} associated to $[\phe]$.

\me We have explicit (and easy to use) techniques to compute the toric degree and toric $r$-tuples (reduced in the torsion case) of $[\phe]$. Furthermore, we can also prove that, in the torsion case, the torsion of $e^{2\pi i [\phe]}$ always divides $m$ (see Proposition 5.1).

\sm As expected, we are able to show that the torsion-free case behaves as the
vector fields case, proving the following analogue of Theorem \rf{ZungIntro} (which 
works even when $\d f_O$ is not diagonalizable).

\sm\thm{Te1.1Intro}{Theorem} {\sl Let~$f$ be a germ of biholomorphism of~$\C^n$ fixing the origin~$O$, such that, denoted by $\lambda=\{\lambda_1, \dots, \lambda_n\}$ the spectrum of~$\d f_O$, the unique $[\phe]\in(\C/\Z)^n$ such that~$\lambda=e^{2\pi i [\phe]}$ is of toric degree~$1\le r\le n$ and torsion-free. Then~$f$ admits a holomorphic Poincar\'e-Dulac normalization if and only if there exists a holomorphic effective action on~$(\C^n,O)$ of a torus of dimension~$r$ commuting with~$f$ and such that the columns of the weight matrix of the action are a~$r$-tuple of toric vectors associated to~$[\phe]$.}

\me The torsion case is more delicate and difficult to deal. First of all, given $[\phe]\in(\C/\Z)^n$ with toric degree $1\le r\le n$ and torsion $\tau\ge2$, and a reduced toric $r$-tuple $\eta^{(1)},\dots, \eta^{(r)}$, we always have
$$ 
\bigcap_{k=2}^r \res_j^+(\eta^{(k)}) \supseteq \res_j([\phe])\supseteq
\bigcap_{k=1}^r \res_j^+(\eta^{(k)}).\tag{trecasiintro}
$$
This suggest a subdivision in several subcases, all realizable (we have examples for all of them) and, surprisingly, having very different behaviours. We have cases similar to the case of germs of vector fields (even if we have torsion!), and cases that are indeed different.
In particular, considering iterates of $f$ to reduce to the torsion-free case hides very interesting phenomena, and it does not allow to see that some torsion cases can be directly studied. 
Moreover, we have explicit (and computable) techniques to decide in which subcase a given $[\phe]\in(\C/\Z)^n$ is.

\me\defin{ImpureTorsionIntro} Let~$[\phe]\in(\C/\Z)^n$ be of toric degree~$1\le r\le n$ and in the torsion case. We say that $[\phe]$ is in the {\it impure torsion case} if, for one (and hence any: see Lemma 7.6) reduced $r$-tuple $\eta^{(1)},\ldots,\eta^{(r)}$ of toric vectors associated to $[\phe]$ we have
$$
\res_j([\phe]) =\bigcap_{k=2}^r\res_j^+(\eta^{(k)}),\tag{torsioneimpura2}
$$
for any $j\in\{1,\dots, n\}$. Otherwise we say that $[\phe]$ is in the {\it pure torsion case}.

\me The impure torsion case is the subcase behaving as the case of germs of vector fields, and in which, again, we do not need $\d f_O$ diagonalizable. In fact, we can prove the following

\sm\thm{Pr4.1Intro}{Theorem} {\sl Let~$f$ be a germ of biholomorphism of~$\C^n$ fixing the origin~$O$, such that, denoted by $\lambda=\{\lambda_1, \dots, \lambda_n\}$ the spectrum of~$\d f_O$, the unique $[\phe]\in(\C/\Z)^n$ such that~$\lambda=e^{2\pi i [\phe]}$ is of toric degree~$1\le r\le n$ and in the impure torsion case. Then it admits a holomorphic Poincar\'e-Dulac normalization if and only if there exists a holomorphic effective action on~$(\C^n,O)$ of a torus of dimension~$r-1$ commuting with~$f$, and such that the columns of the weight matrix of the action are reduced torsion-free toric vectors associated to $[\phe]$.}

\sm The next subcase is

\sm\defin{ReduciblePureTorsionIntro} Let~$[\phe]\in(\C/\Z)^n$ be of toric degree~$1\le r\le n$ and in the pure torsion case. We say that $[\phe]$ {\it can be simplified} if it admits a reduced~$r$-tuple of toric vectors~$\eta^{(1)},\dots, \eta^{(r)}$ such that
$$
\res_j([\phe]) =\bigcap_{k=1}^r\res_j^+(\eta^{(k)}),\tag{torsionepura3}
$$
for all $j=1,\dots, n$. The $r$-tuple $\eta^{(1)},\dots, \eta^{(r)}$ is said a {\it simple reduced} $r$-tuple associated to $[\phe]$.

\sm As we shall see in section 7, condition \rf{torsionepura3} depends on the chosen toric $r$-tuple. However, we have techniques to decide whether $[\phe]$ can be simplified or not.

\sm The case in which $[\phe]$ can be simplified is similar to the case of germs of vector fields, but we have a distinction between the case of $\d f_O$ diagonalizable and $\d f_O$ not diagonalizable, as we see in the following result (for the definition of compatibility see section 3):

\sm\thm{Pr4.2Intro}{Theorem} {\sl Let~$f$ be a germ of biholomorphism of~$\C^n$ fixing the origin~$O$, such that, denoted by $\lambda=\{\lambda_1, \dots, \lambda_n\}$ the spectrum of~$\d f_O$, the unique $[\phe]\in(\C/\Z)^n$ such that~$\lambda=e^{2\pi i [\phe]}$ is of toric degree~$1\le r\le n$ and in the pure torsion case and it can be simplified. Then:{\parindent=30pt
\sm\item{(i)} if~$\d f_O$ is diagonalizable, $f$ admits a holomorphic Poincar\'e-Dulac normalization if and only if there exists a holomorphic effective action on~$(\C^n,O)$ of a torus of dimension~$r$ commuting with~$f$ and such that the columns of the weight matrix $\Theta$ of the action are a simple reduced~$r$-tuple of toric vectors associated to~$[\phe]$;
\sm\item{(ii)} if~$\d f_O$ is not diagonalizable and there exists a simple reduced~$r$-tuple of toric vectors associated to~$[\phe]$ such that its vectors are the columns of a matrix $\Theta$ compatible with~$\d f_O$,~$f$ admits a holomorphic Poincar\'e-Dulac normalization if and only if there exists a holomorphic effective action on~$(\C^n,O)$ of a torus of dimension~$r$ commuting with~$f$ and with weight matrix $\Theta$.
\sm}}

\me The case in which $[\phe]$ cannot be simplified  is the furthest from the case of germs of vectors fields, because we cannot reduce the multiplicative resonances to additive ones. In fact, we only have a sufficient condition for holomorphic normalization.

\sm\thm{PrPureTorsionIntro}{Proposition} {\sl Let~$f$ be a germ of biholomorphism of~$\C^n$ fixing the origin~$O$, such that, denoted by $\lambda=\{\lambda_1, \dots, \lambda_n\}$ the spectrum of~$\d f_O$, the unique $[\phe]\in(\C/\Z)^n$ such that~$\lambda=e^{2\pi i [\phe]}$ is of toric degree $1\le r\le n$ and in the pure torsion case and it cannot be simplified. If there exists a holomorphic effective action on~$(\C^n,O)$ of a torus of dimension~$r$ commuting with~$f$ and such that the columns of the weight matrix of the action are a reduced~$r$-tuple of toric vectors associated to~$[\phe]$, then $f$ admits a holomorphic Poincar\'e-Dulac normalization.}

\me We have then completely understood the relations between torus actions, holomorphic Poin\-ca\-r\'e-Dulac normalizations, and torsion phenomena. We end the article giving an example of techniques to construct torus actions.

\sm\defin{De4.3Intro} Let $1\le m\le n$. A {\it set of $m$ integrable vector fields} of $(\C^n, O)$ is a set~$X_1, \ldots, X_m$ of germs of holomorphic vector fields of~$(\C^n, O)$ singular at the origin, of order $1$ and such that:{\parindent=30pt
\sm\item{(i)}~$X_1, \ldots, X_m$ commute pairwise and are linearly independent;
\sm\item{(ii)} there exist~$n-m$ germs of holomorphic functions~$g_1, \ldots, g_{n-m}$ in~$(\C^n, O)$ which are common first integrals of~$X_1, \ldots, X_m$, and they are functionally independent almost everywhere.
\sm}

\sm\defin{De5.1Intro} A germ of biholomorphism~$f$ of~$(\C^n, O)$ fixing the origin {\it commutes with a set of integrable vector fields} if there exists a positive integer~$1\le m\le n$, such that there exists a set of~$m$ germs of holomorphic integrable vector fields~$X_1, \ldots, X_m$ such that 
$$
\d f(X_j)=X_j\circ f
$$
for each~$j=1, \dots, m$.

\sm A germ of biholomorphism $f$ of $(\C^n,O)$ commutes with a
vector field $X$ according to the previous definition if and only if it commutes with the flow generated by $X$. Then (see also Section 8 for more general
statements):

\sm\thm{Te5.2Intro}{Theorem} {\sl Let~$f$ be a germ of biholomorphism of~$(\C^n, O)$ fixing the origin and commuting with a set of integrable holomorphic vector fields $X_1,\dots, X_m$. Then~$f$ commutes with a holomorphic effective action on~$(\C^n,O)$ of a torus of dimension equal to the toric degree $r$ of~$X_1$ and such that the columns of the weight matrix of the action are a~$r$-tuple of toric vectors associated to~$X_1$.}

\me If the semi-simple part of the linear term of $X_1$ has the same resonances of the eigenvalues of $\d f_O$, then we get a Poincar\'e-Dulac holomorphic normalization. Moreover, the commutation condition implies that $f$ preserves the foliation generated by $X_1, \dots, X_m$; hence this is a condition similar to condition A of [Brj] for the convergence of a normalization in the case of vector fields. Finally, commuting with a torus action implies that $f$ preserves the orbits of the action, hence a condition of the kind of Definition \rf{De5.1Intro} is close to be necessary and sufficient for the existence of a commuting torus action. 

\me The structure of this paper is as follows. 

In section $2$ we shall recall some basic facts about linear torus actions and we shall fix some notations.
In section $3$ we shall describe the relations between the existence of torus actions with certain properties and the possibility to conjugate a given germ of biholomorphism to another one of a particular form, and we shall prove Theorem \rf{Toro} and Theorem \rf{Linearizzazione}.
In section $4$ we shall study the toric degree and the toric $r$-tuples associated to the eigenvalues of the differential~$\d f_O$ of a germ of biholomorphism of $\C^n$ fixing the origin, and the weight matrices of torus actions.  
In section $5$ we shall study the notion of torsion and we shall prove Theorem \rf{PrNoImpureTorsionFreeIntro}.
In section $6$ we shall study the relations between the resonances of the eigenvalues of the differential~$\d f_O$ of a germ of biholomorphism of $\C^n$ fixing the origin and the additive resonances of an associated toric $r$-tuple, in the torsion-free case, and we shall prove Theorem \rf{Te1.1Intro}.
In section $7$ we shall study the relations between the resonances of the eigenvalues of the differential~$\d f_O$ of a germ of biholomorphism of $\C^n$ fixing the origin and the additive resonances of an associated toric $r$-tuple, in the torsion case, and we shall prove Theorem \rf{Pr4.1Intro}, Theorem~\rf{ReduciblePureTorsionIntro} and Proposition \rf{PrPureTorsionIntro}.
In the last section we shall give some geometric conditions to construct the torus actions we need, and we shall prove Theorem \rf{Te5.2Intro} and other similar results.

\me{\bf Acknowledgments.} I would like to thank Marco Abate for helpful comments on a draft of this work, and Jean \'Ecalle for the useful conversations on the last part of section $7$. Part of this work was done during a visit to the D\'epartement de Math\'ematiques de la Facult\'e des Sciences d'Orsay, Universit\'e Paris-Sud 11, and I would also like to thank the D\'epartement d'Orsay for its hospitality.

%%%%%%%%%%%%%%%%%%%%%%%%%%%%%%%%%%%%%%%%%%%%%%%%%%%%%%%%%%%%%%%%%%%%%%%%%%%%%%%%
\sect Preliminaries
%%%%%%%%%%%%%%%%%%%%%%%%%%%%%%%%%%%%%%%%%%%%%%%%%%%%%%%%%%%%%%%%%%%%%%%%%%%%%%%%

%{\sf qui metto le definizioni varie di risonanze additive e moltiplicative, sistemare la terminologia sulle azioni di tori etc $|Q|=\sum_{h=1}^n q_h\ge 2$ $\langle Q,\theta\rangle =\sum_{h=1}^n q_h\theta_h $ say that $[\cdot]\colon
%\C^n\to(\C/\Z)^n$ is the standard projection.}

\sm Let $A\colon\T^r\times M\to M$ be a torus action on a complex manifold~$M$,
with a fixed point $p_0\in M$ (that is $A(x,p_0)=A_x(p_0)=p_0$ for all
$x\in\T^r$). The differential $\d(A_x)_{p_0}\colon T_{p_0}M\to T_{p_0}M$
is then well-defined, and thus we have a linear torus action on~$T_{p_0}M$.
A linear torus action can be thought of as a Lie group homomorphism
$A\colon\T^r\to\Aut(T_{p_0}M)$, that is as a representation of $\T^r$
on the vector space $V=T_{p_0}M$. 

Characters and weights of $\T^r$ are well known. All characters of
$\T^1=\S^1=\R/\Z$ are of the form 
$$
\chi_\theta(x)=\exp(2\pi i x\theta)
$$
with $\theta\in\Z$; hence the character group of~$\T^1$ is isomorphic to~$\Z$. Since $\T^r=\T^1\times\cdots\times\T^1$,
the characters of~$\T^r$ are obtained multiplying characters of~$\T^1$,
that is they are of the form
$$
\chi_\theta(x)=\exp\left(2\pi i\sum_{k=1}^r x_k\theta^k\right),
$$ 
with $\theta=(\theta^1,\ldots,\theta^r)\in(\Z^r)^*$, where the ${}^*$ 
denotes the dual. In particular, $\theta$ should be thought of as
a row vector. 
The weights of $\T^r$ are then the differential of the characters
computed at the identity element, and thus are given by
$$
w_\theta(v)=2\pi i\sum_{k=1}^r v_k\theta^k
$$
with $\theta\in(\Z^r)^*$ and $v\in\R^r$. 
If we write $\theta_j=(\theta_j^1,\ldots,\theta_j^r)\in(\Z^r)^*$, then the
matrix representation of the linear action~$A$ in the eigenvector basis is given by
$$
A(x)=\hbox{diag}\left(\chi_{\theta_j}(x)\right)=\hbox{diag}\left(\exp\left(
2\pi i\sum_{k=1}^r x_k\theta^k_j\right)\right).
$$
We have then associated to our torus action a matrix $\Theta=(\theta^k_j)\in M_{n\times r}(\Z)$, whose columns do not depend on the particular coordinates chosen to express the torus action, but can be uniquely (up to order) recovered by the action itself.

\sm\defin{MatricePesi} The matrix $\Theta$ just defined is called the {\it weight matrix} of the torus action.

\sm\defin{De2.0} Let~$\theta\in\C^n$ and let~$j\in\{1, \dots, n\}$. We say that a multi-index~$Q\in\N^n$, with~$|Q|=\sum_{h=1}^n q_h\ge 2$, gives an {\it additive resonance relation for~$\theta$ relative to the~$j$-th coordinate} if
$$
\langle Q,\theta\rangle =\sum_{h=1}^n q_h\theta_h = \theta_j
$$ 
and we put
$$
{\rm Res}_j^+(\theta)=\{Q\in\N^n\mid |Q|\ge 2,\langle Q,\theta\rangle = \theta_j\}.
$$
%We say that a multi-index~$Q\in\N^n$, with~$|Q|\ge 2$, gives a {\it multiplicative resonance relation for~$\theta$ relative to the~$j$-th coordinate} if
%$$
%\langle Q,\theta\rangle - \theta_j\in\Z
%$$ 
%and we put
%$$
%{\rm Res}_j(\theta)=\{Q\in\N^n\mid |Q|\ge 2,\langle Q,\theta\rangle - \theta_j\in\Z\}.
%$$
Let~$\lambda\in(\C^*)^n$ and let~$j\in\{1, \dots, n\}$. We say that a multi-index~$Q\in\N^n$, with~$|Q|\ge 2$, gives a {\it multiplicative resonance relation for~$\lambda$ relative to the~$j$-th coordinate} if
$$
\lambda^Q =\lambda_1^{q_1}\cdots\lambda_n^{q_n} = \lambda_j
$$ 
and we put
$$
{\rm Res}_j(\lambda)=\{Q\in\N^n\mid |Q|\ge 2,\lambda^Q = \lambda_j\}.
$$ 

\sm\thm{ReMolt1}{Remark} Given $[\phe]\in(\C/\Z)^n$, where $[\,\cdot\,]\colon
\C^n\to(\C/\Z)^n$ is the standard projection, the set
$$
\{Q\in\N^n\mid |Q|\ge 2,\langle Q,\phe\rangle - \phe_j\in\Z\}.
$$
does not depend on the specific representative
$\phe\in\C^n$ but only on the class~$[\phe]$, and so it is well defined the set 
$\res_j([\phe])$ as
$$
\res_j([\phe])=\{Q\in\N^n\mid |Q|\ge 2,\left[\langle Q,\phe\rangle - \phe_j\right]=[0]\}.
$$

\sm\thm{ReMolt2}{Remark} Notice that given~$\lambda\in(\C^*)^n$, we can always find a unique~$[\phe]\in(\C/\Z)^n$ such that~$\lambda=e^{2\pi i [\phe]}$, i.e.,~$\lambda_j = e^{2\pi i [\phe_j]}$ for every~$j=1, \dots, n$. Then $\res_j(\lambda)=\res_j([\phe])$, 
thus justifying the definitions and the terminology. 

%%%%%%%%%%%%%%%%%%%%%%%%%%%%%%%%%%%%%%%%%%%%%%%%%%%%%%%%%%%%%%%%%%%%%%%%%%%%%%%%
\sect Torus Actions and Normal Forms of germs of biholomorphisms
%%%%%%%%%%%%%%%%%%%%%%%%%%%%%%%%%%%%%%%%%%%%%%%%%%%%%%%%%%%%%%%%%%%%%%%%%%%%%%%%

\sm In this section we shall describe the relations between the existence of torus actions with certain properties and the possibility to conjugate a given germ of biholomorphism to another of a particular form.

\sm\defin{De2.1} Let~$\theta^{(1)},\dots, \theta^{(r)}\in\Z^n$. We say that a monomial~$z^Qe_j$, with~$Q\in\N^n$,~$|Q|\ge 1$ and~$j\in\{1,\dots, n\}$, is{\it~$\Theta$-resonant}, where $\Theta$ is the~$n\times r$ matrix whose columns are~$\theta^{(1)},\dots, \theta^{(r)}$, if
$$
\langle Q,\theta^{(k)}\rangle = \theta^{(k)}_j
$$ 
for every~$k=1, \dots, r$, that is if $\theta^{(k)}_h=\theta^{(k)}_j$, for $Q=e_h$, or
$$
Q\in\ca R_j(\Theta)= \bigcap_{k=1}^r\res_j^+(\theta^{(k)}),\tag{eqzero}
$$
for $|Q|\ge 2$. We say that $\Theta$ has {\it no resonances} if $\ca R_j(\Theta)= \void$  for every~$j=1,\dots, n$.

\sm\defin{De2.2} Let~$\theta^{(1)},\dots, \theta^{(r)}\in\Z^n$ and let~$T$ be a linear map of $\C^n$. We say that the matrix $\Theta$, with columns~$\theta^{(1)},\dots, \theta^{(r)}$, is {\it compatible with} $T$ if and only if we can write $T$ in Jordan form with all monomials
$\Theta$-resonant.

\sm\thm{Te2.1}{Theorem} {\sl Let $f$ be a germ of biholomorphism of $\C^n$ fixing the origin~$O$. Then $f$ commutes with a holomorphic effective action on $(\C^n,O)$ of a torus of dimension $1\le r\le n$ with weight matrix $\Theta\in M_{n\times r}(\Z)$ if and only there exists a local holomorphic change of coordinates conjugating $f$ to a germ with linear part in Jordan normal form and containing only $\Theta$-resonant monomials.}

\sm\proof Let us suppose that the linear part of~$f$ is in Jordan normal form and~$f$ contains only~$\Theta$-resonant monomials. Then we claim that~$f$ commutes with the linear effective torus action
$$
A\colon\T^r\times (\C^n, O)\to(\C^n, O),
$$ 
defined by
$$
A(x, z) = \diag\left(e^{2\pi i \sum_{k=1}^r x_k \theta^{(k)}_j}\right)z.
$$ 
In fact in these hypotheses the~$j$-th coordinate of~$f$ is
$$
\lambda_j z_j + \eps_j z_{j-1} + \sum_{|Q|\ge 2\atop Q\in\ca R_j(\Theta)} f_{Q,j} z^Q
$$
where~$\eps_j\in\{0,1\}$ can be different from~$0$ only if~$\lambda_j=\lambda_{j-1}$, the set~$\ca R_j(\Theta)$ is defined in \rf{eqzero} and the assumption that $\eps_j z_{j-1}e_j$ is $\Theta$-resonant implies $\theta^{(k)}_{j-1}=\theta^{(k)}_j$ for $k=1,\ldots,r$ if~$\eps_j\ne 0$. Then for every~$x\in\T^r$ we have
$$
\eqalign{
f_j(A(x, z)) &= \lambda_j e^{2\pi i\sum_{k=1}^r x_k \theta^{(k)}_j}z_j + \eps_j e^{2\pi i\sum_{k=1}^r x_k \theta^{(k)}_{j-1}}z_{j-1} + \!\!\!\sum_{|Q|\ge 2\atop Q\in\ca R_j(\Theta)} f_{Q,j}e^{2\pi i\sum_{k=1}^r x_k\langle Q,\theta^{(k)}\rangle} z^Q \cr 
							&=\lambda_j e^{2\pi i\sum_{k=1}^r x_k \theta^{(k)}_j}z_j + \eps_j e^{2\pi i\sum_{k=1}^r x_k \theta^{(k)}_{j}}z_{j-1} +\!\!\! \sum_{|Q|\ge 2\atop Q\in\ca R_j(\Theta)} f_{Q,j}e^{2\pi i\sum_{k=1}^r x_k\theta^{(k)}_j} z^Q \cr
							&=e^{2\pi i\sum_{k=1}^r x_k \theta^{(k)}_j}\left(\lambda_j z_j + \eps_j z_{j-1} + \sum_{|Q|\ge 2\atop Q\in\ca R_j(\Theta)} f_{Q,j} z^Q\right)\cr
							&=e^{2\pi i\sum_{k=1}^r x_k \theta^{(k)}_j}\left(f_j(z)\right)\cr
							&=A(x, f(z))_j. 
}
$$

Conversely, let us suppose that~$f$ commutes with a holomorphic effective action on~$(\C^n,O)$ of a torus of dimension~$1\le r\le n$ with weight matrix~$\Theta$. Then, by B\"ochner linearization theorem [B], there exists a tangent to the identity holomorphic change of variables~$\psi$ linearizing the torus action. Furthermore, up to a linear change of coordinates we can assume that in the new coordinates the action is given by 
$$
A(x, z) = \diag\left(e^{2\pi i \sum_{k=1}^r x_k \theta^{(k)}_j}\right)z,
$$ 
and that $f$ (still commuting with the torus action) has
linear part in Jordan normal form compatible with $\Theta$, and thus its~$j$-th coordinate is
$$
\lambda_j z_j + \eps_j z_{j-1} + \sum_{|Q|\ge 2} f_{Q,j} z^Q
$$
where~$\eps_j\in\{0,1\}$ can be different from~$0$ only if $\lambda_{j-1}=\lambda_j$ and $\theta_{j-1}=\theta_j$. For every~$x\in\T^r$, we have
$$
\eqalign{
f_j( A(x, z)) &= \lambda_j e^{2\pi i\sum_{k=1}^r x_k \theta^{(k)}_j}z_j + \eps_j e^{2\pi i\sum_{k=1}^r x_k \theta^{(k)}_{j-1}}z_{j-1} + \sum_{|Q|\ge 2} f_{Q,j}e^{2\pi i\sum_{k=1}^r x_k\langle Q,\theta^{(k)}\rangle} z^Q,}
$$ %\cr 
							%&=\lambda_j e^{2\pi i\sum_{k=1}^r x_k \theta^{(k)}_j}z_j + \eps_j e^{2\pi i\sum_{k=1}^r x_k \theta^{(k)}_{j}}z_{j-1} + \sum_{|Q|\ge 2} f_{Q,j}e^{2\pi i\sum_{k=1}^r x_k\langle Q,\theta^{(k)}\rangle} z^Q ,}$$
and
$$
A(x, f(z))_j=e^{2\pi i\sum_{k=1}^r x_k \theta^{(k)}_j}\left(\lambda_j z_j + \eps_j z_{j-1} + \sum_{|Q|\ge 2} f_{Q,j} z^Q\right).
$$
Then~$f_j( A(x, z))=A(x, f(z))_j$ if and only if
$$
f_{Q,j}\left(e^{2\pi i\sum_{k=1}^r x_k\langle Q,\theta^{(k)}\rangle}-e^{2\pi i\sum_{k=1}^r x_k \theta^{(k)}_j} \right)=0
$$ 
for every~$x\in\T^r$,~$j=1,\cdots,n$,~$Q\in\N^n$ with~$|Q|\ge 2$, i.e.,~$f_{Q,j}$ can be non-zero only when
$$
\sum_{k=1}^r x_k(\langle Q, \theta^{(k)}\rangle-\theta^{(k)}_j)\in\Z \quad\forall x\in\T^r,
$$
which is equivalent to
$$
\langle Q, \theta^{(k)}\rangle-\theta^{(k)}_j=0
$$
for every~$k=1,\dots, r$, meaning that~$f$ contains only~$\Theta$-resonant monomials. \qed

\sm As a consequence of this result we have

\sm\thm{Co2.1}{Corollary} {\sl Let~$f$ be a germ of biholomorphism of~$\C^n$ fixing the origin~$O$. Then~$f$ is holomorphically linearizable if and only if it commutes with a holomorphic effective action on~$(\C^n,O)$ of a torus of dimension~$1\le r\le n$ with weight matrix~$\Theta$ having no resonances.}

\sm\proof If $f$ is linear and in Jordan normal form, then it
commutes with any linear action of~$\T^1$ with compatible 
weight matrix~$\Theta$; so it suffices to choose $\Theta$ with $\ca R_1(\Theta)=\ldots=\ca R_n(\Theta)=\void$.% \`e ovvio che ne esiste almeno uno infatti basta sceglierne uno i cui multi-indici di risonanza non siano ammissibili.

Conversely, if~$f$ commutes with a holomorphic effective action on~$(\C^n,O)$ of a torus of dimension~$1\le r\le n$ with weight matrix~$\Theta$, then, by the previous result, $\Theta$ is compatible with the linear part of $f$ and there exists a local holomorphic change of coordinates such that~$f$ is conjugated to a germ with the same linear part and containing only~$\Theta$-resonant monomials. But each~$\ca R_j(\Theta)$ is empty; hence there are no~$\Theta$-resonant monomials of degree at least $2$, and thus~$f$ is holomorphically linearizable. \qed

%%%%%%%%%%%%%%%%%%%%%%%%%%%%%%%%%%%%%%%%%%%%%%%%%%%%%%%%%%%%%%%%%%%%%%%%%%%%%%%%
\sect Toric degree 
%%%%%%%%%%%%%%%%%%%%%%%%%%%%%%%%%%%%%%%%%%%%%%%%%%%%%%%%%%%%%%%%%%%%%%%%%%%%%%%%

\sm We want to study the relations between the resonances of the eigenvalues of the differential~$\d f_O$ of a germ of biholomorphism of $\C^n$ fixing the origin, and the weight matrices of torus actions to understand in which cases Theorem \rf{Te2.1} gives us a Poincar\'e-Dulac holomorphic normalization. Thanks to Remark \rf{ReMolt2} we have to deal with vectors of $(\C/\Z)^n$. A concept that turns out to be crucial for this study is that of {\it toric degree}.

\sm\defin{De1.1} Let~$[\phe]=([\phe_1],\dots, [\phe_n])\in(\C/\Z)^n$. The {\it toric degree of~$[\phe]$} is the minimum~$r\in\N$ such that there exist~$\alpha_1, \dots, \alpha_r\in\C$ and~$\theta^{(1)},\dots, \theta^{(r)}\in\Z^n$ such that
$$
[\phe]= \left[\sum_{k=1}^r \alpha_k \theta^{(k)}\right].\tag{Eq1.1}
$$ 
The~$r$-tuple~$\theta^{(1)},\dots, \theta^{(r)}$ is called a {\it $r$-tuple of toric vectors associated to~$[\phe]$}, and the numbers~$\alpha_1, \dots, \alpha_n\in\C$ are {\it toric coefficients} of the toric $r$-tuple.

\sm\thm{ReToric}{Remark} Note that the toric degree is necessarily at most $n$, since
$$
[\phe]= \left[\sum_{k=1}^n \phe_k e_k\right].
$$ 

\sm We did not say {\it the} toric coefficients because we have the following result.

\sm\thm{Re1.4}{Lemma} {\sl Let~$[\phe]\in(\C/\Z)^n$ be of toric degree $1\le r\le n$ and let~$\theta^{(1)},\dots, \theta^{(r)}$ be a~$r$-tuple of toric vectors associated to~$[\phe]$ with toric coefficients~$\alpha_1,\dots, \alpha_r\in\C$. Then~$\beta_1, \dots, \beta_r\in\C$ satisfy
$$
[\phe]= \left[\sum_{k=1}^r \beta_k \theta^{(k)}\right]
$$ 
if and only if
$$
\Theta\pmatrix{\alpha_1-\beta_1\cr
				\vdots\cr
				\alpha_r-\beta_r} \in\Z^n
$$
where $\Theta$ is the~$n\times r$ matrix whose columns are~$\theta^{(1)},\dots, \theta^{(r)}$.}

\sm\proof We have
$$
\left[\sum_{k=1}^r \alpha_k \theta^{(k)}\right] = \left[\sum_{k=1}^r \beta_k \theta^{(k)}\right]
$$
if and only if 
$$
\sum_{k=1}^r \alpha_k \theta^{(k)}-\sum_{k=1}^r \beta_k \theta^{(k)}\in\Z^n,
$$
that is
$$
\Theta\pmatrix{\alpha_1-\beta_1\cr
				\vdots\cr
				\alpha_r-\beta_r} \in\Z^n,
$$
which is the assertion. \qed

\sm Thanks to Remark \rf{ReMolt2} the following definition makes sense.

\sm\defin{De1.2} Let~$f$ be a germ of biholomorphism of~$\C^n$ fixing the origin and denote by $\lambda=\{\lambda_1, \dots, \lambda_n\}$ the spectrum of~$\d f_O$. The {\it toric degree of~$f$} is the toric degree of the unique~$[\phe]\in(\C/\Z)^n$ such that~$\lambda=e^{2\pi i [\phe]}$ .

\sm Toric $r$-tuples and toric coefficients have to satisfy certain arithmetic properties, as the following result shows. 

\sm\thm{Le1.0}{Lemma} {\sl Let~$[\phe]\in(\C/\Z)^n$ be of toric degree~$1\le r\le n$ and let~$\theta^{(1)},\dots, \theta^{(r)}$ be a~$r$-tuple of toric vectors associated to~$[\phe]$ with toric coefficients~$\alpha_1,\dots, \alpha_r\in\C$. Then:{\parindent=30pt
\sm\item{(i)}~$\alpha_1, \dots, \alpha_r$ is a set of rationally independent complex numbers;
\sm\item{(ii)} every~$r$-tuple of toric vectors associated to~$[\phe]$ is a set of~$\Q$-linearly independent vectors.
\sm}}

\sm\proof (i) Let us suppose by contradiction that~$\alpha_1, \dots, \alpha_r\in\C$ are rationally dependent. Then there exists~$(c_1,\dots, c_r)\in\Z^r\setminus\{O\}$ such that
$$
c_1\alpha_1 +\cdots+ c_r\alpha_r =0.
$$
Up to reordering we may assume~$c_1\ne 0$. Then
$$
\alpha_1= -{1\over c_1} (c_2\alpha_2 +\cdots+ c_r\alpha_r),
$$
and hence
$$
\eqalign{
[\phe] 	&=\left[\sum_{k=1}^r \alpha_k \theta^{(k)}\right]\cr
			  &= \left[-{1\over c_1} (c_2\alpha_2 +\cdots+ c_r\alpha_r)\theta^{(1)} + \alpha_2\theta^{(2)} + \cdots+ \alpha_r\theta^{(r)}\right]\cr
			  &= \left[{\alpha_2\over c_1}(c_1\theta^{(2)}-c_2\theta^{(1)})+\cdots+ {\alpha_r\over c_1}(c_1\theta^{(r)}-c_r\theta^{(1)})\right],}
$$ 
and this contradicts the definition of toric degree. 

\sm (ii) The proof is analogous to the previous one. \qed

\sm\thm{Re1.0}{Remark} Given~$[\phe]\in(\C/\Z)^n$, of toric degree~$1\le r\le n$, if~$\theta^{(1)},\dots, \theta^{(r)}$ is a~$r$-tuple of toric vectors associated to~$[\phe]$, the~$n\times r$ matrix~$\Theta$ whose columns are~$\theta^{(1)},\dots, \theta^{(r)}$, has maximal rank~$r$. 

\sm\thm{ReCoprimalita}{Remark} Note that, if~$[\phe]\in(\C/\Z)^n$ has toric degree~$1\le r\le n$, and~$\theta^{(1)},\dots, \theta^{(r)}$ is a~$r$-tuple of toric vectors associated to~$[\phe]$, up to change the toric coefficients $\alpha_1, \dots, \alpha_r$, we can always assume $\theta^{(k)}_1, \dots, \theta^{(k)}_n$ coprime for each $1\le k\le r$. In fact, if $d_k\in\Z$ is the greatest common divisor of $\theta^{(k)}_1, \dots, \theta^{(k)}_n$, then 
$$
[\phe]= \left[\sum_{k=1}^r \alpha_k\theta^{(k)}\right] = \left[\sum_{k=1}^r d_k\alpha_k\widetilde\theta^{(k)}\right],
$$
where
$$
\widetilde\theta^{(k)}=\pmatrix{\theta^{(k)}_1/d_k\cr
								\vdots\cr
								\theta^{(k)}_n/d_k}
$$
for $k=1,\dots, r$.

\sm\thm{Re1.5}{Remark} {\sl Given~$[\phe]\in(\C/\Z)^n$, of toric degree~$1\le r\le n$, the $r$-tuple of toric vectors associated to $[\phe]$ is not necessarily unique.} Let us consider, for example
$$
[\phe]= \left[\matrix{{3\sqrt{2} + 4i} \cr
		   2\sqrt{2} + 6i \cr
		   -\sqrt{2} + 2i }\right].
$$
The toric degree cannot be~$1$, since it is immediate to verify that~$\phe$ cannot be written as the product of a complex number times an integer vector.
The toric degree is in fact~$2$, since we have
$$
[\phe]= \left[\sqrt{2}\pmatrix{3 \cr
		   2 \cr
		   -1 } + 2i\pmatrix{2 \cr
		   3 \cr
		   1 }
		   \right].
$$
However we can also write $[\phe]$ as
$$
[\phe]= \left[{-3\sqrt{2}+ 16i\over 6}\pmatrix{ 0\cr
		   1 \cr
		   1 } + {3\sqrt{2} + 4i \over 6}\pmatrix{6 \cr
		   5 \cr
		   -1 }
		   \right].
$$
Note that, in both cases, the toric coefficients are rationally independent with $1$.

\sm\thm{Ex1.3}{Example} The vector of $(\C/\Z)^2$
$$
[\phe]= \left[\matrix{{(1+ 6\sqrt{2})/6} \cr
					{(1-2\sqrt{2})/ 2}}\right],
$$
has toric degree $2$, since we have
$$
[\phe]= \left[{1+ 6\sqrt{2}\over 6} \pmatrix{1 \cr
		   0} + {1-2\sqrt{2}\over 2}\pmatrix{0 \cr
		   1}
		   \right],
$$
and it is not difficult to verify that it cannot have toric degree $1$. We can also write $[\phe]$ as
$$
[\phe]= \left[{1\over 6}\pmatrix{1 \cr
		   3 } + \sqrt{2}\pmatrix{1 \cr
		   -1}
		   \right].
$$
Note that this time, in both cases, the toric coefficients are rationally dependent with $1$.

\sm We shall prove that, given~$[\phe]\in(\C/\Z)^n$ of toric degree~$1\le r\le n$, even when the $r$-tuple of toric vectors associated to $[\phe]$ is not unique, we can always say whether the toric coefficients are rationally independent with $1$ or not, so this will be an intrinsic property of the vector.
Before proving this, we shall need the following result that gives us a way to find a more useful toric $r$-tuple when the toric coefficients are rationally dependent with $1$.

\sm\thm{ReGrado1}{Remark} Note that~$\alpha\in\C$ is rationally dependent with $1$ if and only if it belongs to $\Q$.

\sm\thm{Pr1.0}{Lemma} {\sl Let~$[\phe]\in(\C/\Z)^n$ be of toric degree~$1\le r\le n$, and let~$\theta^{(1)},\dots, \theta^{(r)}$ be a~$r$-tuple of toric vectors associated to~$[\phe]$ with toric coefficients~$\alpha_1,\dots, \alpha_r\in\C$ rationally dependent with $1$. Then there exists a~$r$-tuple of toric vectors~$\eta^{(1)},\dots, \eta^{(r)}$ associated to~$[\phe]$ with toric coefficients~$\beta_1, \dots, \beta_r\in\C$ such that $\beta_1=1/m$ with~$m\in\N\setminus\{0,1\}$ and~$m, \eta^{(1)}_1, \dots, \eta^{(1)}_n$ coprime. Moreover $\beta_2, \dots, \beta_r$ are rationally independent with $1$.}

\sm\proof If~$r=1$, then~$\alpha$ is rationally dependent with $1$ if and only if it belongs to~$\Q$, i.e.,
$$
[\phe] = \left[{p\over q}\theta\right]
$$
where we may assume without loss of generality~$p$ and~$q$ coprime and~$q, \theta_1,\dots, \theta_n$ coprime. Then
$$
[\phe] = \left[{1\over q} \eta\right]
$$
where~$\eta = p\cdot \theta\in\Z^n$ and we are done. 

Let us suppose now~$r\ge2$. Since~$\alpha_1,\dots,\alpha_r$ are (rationally independent and) rationally dependent with $1$, we can consider the minimum positive integer~$m_0\in\N\setminus\{0\}$ so that there exists~$(m_1, \dots, m_r)\in\Z^r\setminus\{O\}$ such that 
$$
m_1\alpha_1 + \cdots + m_r\alpha_r = m_0.
$$
Thanks to the minimality of $m_0$, we have that $m_1, \dots, m_r, m_0$ are coprime. 
Up to reordering we may assume~$m_1\ne 0$. Then
$$
\eqalign{
\alpha_1&= {m_0\over m_1}-\left({m_2\over m_1}\alpha_2 +\cdots+ {m_r\over m_1}\alpha_r\right)\cr
&={m'_0\over m'_1}-\left({m_2\over m_1}\alpha_2 +\cdots+ {m_r\over m_1}\alpha_r\right),}
$$
where~${m_0\over m_1} = {m'_0\over m_1'}$ with~$(m'_0, m_1')=1$ and~$m_1'\in\N\setminus\{0,1\}$. Let~$d$ be the greatest common divisor of~$m'_1$ and the components of~$\theta^{(1)}$, and consider
$$
\widetilde \theta^{(1)} = {1\over d} \theta^{(1)}, \quad \widetilde m_1 = {m'_1\over d}.
$$ 
Hence
$$
\eqalign{
[\phe] &= \left[{m'_0\over m'_1} \theta^{(1)} + \sum_{k=2}^r{\alpha_k\over m_1} (m_1\theta^{(k)} - m_k\theta^{(1)})\right]\cr
				 &= \left[{m'_0\over \widetilde m_1} \widetilde\theta^{(1)} + \sum_{k=2}^r {\alpha_k\over m_1} (m_1\theta^{(k)} - m_k\theta^{(1)})\right]\cr
				&= \left[{1\over \widetilde m_1} m'_0 \widetilde\theta^{(1)} + \sum_{k=2}^r {\alpha_k\over m_1} (m_1\theta^{(k)} - m_k\theta^{(1)})\right]\cr
			  &= \left[\sum_{k=1}^r \beta_k \eta^{(k)}\right],}
$$ 
where
$$
\beta_1= {1\over \widetilde m_1}, \beta_2 = {\alpha_2\over m_1},\dots, \beta_r = {\alpha_r\over m_1},
$$ 
and
$$
\eta^{(1)} = m'_0 \widetilde\theta^{(1)},\, \eta^{(2)}= m_1\theta^{(2)} - m_2\theta^{(1)},\, \dots,\, \eta^{(r)}= m_1\theta^{(r)} - m_r\theta^{(1)}.
$$
Notice that $\widetilde m_1$ is necessarily greater than $1$, because otherwise the toric degree of $[\phe]$ would be less than $r$.

Now, if~$\beta_2, \dots, \beta_r$ were rationally dependent with $1$, then we would have $k_2,\dots,k_r\in\Z\setminus\{0\}$ such that
$$
k_2\beta_2 + \cdots + k_r\beta_r = k\in\Z\setminus\{0\},
$$
then
$$
- k \,\widetilde m_1 \cdot {1\over \widetilde m_1} + k_2\beta_2 + \cdots + k_r\beta_r = 0,
$$
%by the same argument, we could find~$\gamma_2,\dots, \gamma_r\in\C$ and~$\xi^{(2)}, \dots, \xi^{(r)}\in\Z^n$, with~$\gamma_2= 1/n_2$,~$n_2\in \N\setminus\{0,1\}$ and~$n_2, \xi^{(2)}_1, \dots, \xi^{(2)}_n$ coprime, and
%$$
%\left[\sum_{k=2}^r \beta_k \eta^{(k)}\right] =\left[\sum_{k=2}^r \gamma_k \xi^{(k)}\right],
%$$
%then
%$$
%[\phe] = \left[{1\over m_1'} \eta^{(1)} + {1\over n_2} \xi^{(2)} + \sum_{k=3}^r \gamma_k \xi^{(k)}\right]
%$$
%but
%$$
%\widetilde m_1\cdot {1\over \widetilde m_1} - n_2\cdot {1\over n_2} + 0\cdot \gamma_3 +\cdots + 0\cdot\gamma_r = 0
%$$
contradicting Lemma \rf{Le1.0}. This concludes the proof. \qed

\sm\defin{De4.2} Let~$[\phe]=([\phe_1],\dots, [\phe_n])\in(\C/\Z)^n$ be of toric degree~$1\le r\le n$. We say that a~$r$-tuple~$\eta^{(1)},\dots, \eta^{(r)}$ of toric vectors associated to~$[\phe]$ with rationally dependent with~$1$ toric coefficients~$\beta_1, \dots, \beta_r$ is {\it reduced} if $\beta_1=1/m$ with~$m\in\N\setminus\{0,1\}$ and~$m, \eta^{(1)}_1, \dots, \eta^{(1)}_n$ coprime. In this case the toric vectors $\eta^{(2)},\ldots,\eta^{(r)}$ are called {\it reduced torsion-free toric vectors} associated to $[\phe]$.

\sm Now we can prove that the rational independence with $1$ of the coefficients of toric $r$-tuples associated to a given vector~$[\phe]\in(\C/\Z)^n$ of toric degree~$1\le r\le n$ is an intrinsic property of~$[\phe]$.

\sm\thm{Le1.3}{Proposition} {\sl Let~$[\phe]\in(\C/\Z)^n$ be of toric degree~$1\le r\le n$, and let~$\theta^{(1)},\dots, \theta^{(r)}$ be a~$r$-tuple of toric vectors associated to~$[\phe]$, with toric coefficients~$\alpha_1,\dots, \alpha_r\in\C$ rationally independent with $1$. Then every other~$r$-tuple of toric vectors associated to~$[\phe]$ has toric coefficients rationally independent with $1$.}

\sm\proof Let us assume by contradiction that there exists a~$r$-tuple~$\eta^{(1)},\dots, \eta^{(r)}$ of toric vectors associated to~$[\phe]$ with toric coefficients~$\beta_1,\dots, \beta_r$ rationally dependent with $1$. Thanks to Lemma \rf{Pr1.0}, we may assume without loss of generality~$\beta_1 = 1/m$ with~$m\in\N\setminus\{0,1\}$ and~$m,\eta^{(1)}_1,\dots, \eta^{(1)}_n$ coprime. Let~$N$ be the matrix with columns~$\eta^{(1)},\dots, \eta^{(r)}$, and let~$\Theta$ be the matrix with columns~$\theta^{(1)},\dots, \theta^{(r)}$. We have
$$
[\phe]= \left[N\cdot\pmatrix{\beta_1 \cr
		   \vdots \cr
		   \beta_r}\right] = \left[\Theta\cdot\pmatrix{\alpha_1 \cr
		   \vdots \cr
		   \alpha_r}\right],
$$
that is, there exists an integer vector~${\bf k}\in\Z^n$ such that
$$
N\cdot\pmatrix{\beta_1 \cr
		   \vdots \cr
		   \beta_r} = \Theta\cdot\pmatrix{\alpha_1 \cr
		   \vdots \cr
		   \alpha_r} +{\bf k}.
$$		   
Since~$N$ has maximal rank~$r$, the linear map~$N\colon\Q^r\to\Q^n$ is injective and, for every~$U\subseteq\Q^n$ such that~$\Q^n= \Im(N) \oplus U$, there is a linear map~$L_U\colon\Q^n\to\Q^r$ such that~$\ker(L_U) = U$ and~$L_U N = \id$; hence there is a linear map~$\widetilde L_U\colon\Z^n\to\Z^r$ that~$\widetilde L_U N = h\id$, with~$h\in\Z\setminus\{0\}$.
Then
$$
h\pmatrix{\beta_1 \cr
		   \vdots \cr
		   \beta_r} = \widetilde L_U\Theta\cdot\pmatrix{\alpha_1 \cr
		   \vdots \cr
		   \alpha_r}+\widetilde L_U{\bf k}.
$$
%If~$m$ does not divide~$h$, we are done. In fact~$\widetilde L_U\Theta$ is a matrix with integer coefficients, which implies that~$h\beta_1$ is a linear combination with integer coefficients of~$\alpha_1, \dots, \alpha_r$ and, since~$\beta_1=1/m$, this means that~$\alpha_1, \dots, \alpha_r$ are rationally dependent contradicting the hypothesis.
Moreover, we can choose $U$ so that the first row of~$\widetilde L_U\Theta$ is not identically zero. In fact, the first row of~$\widetilde L_U\Theta$ is identically zero if and only if the first vector~$e_1$ of the standard basis belongs to~$\ker(\Theta^T\widetilde L_U^T)$, and hence it is orthogonal to~$\Im(\widetilde L_U\Theta)$, because for any~$u\in\Q^r$ we have
$$
{0= \la u, \Theta^T\widetilde L_U^T e_1\ra
			= \la\Theta u, \widetilde L_U^T e_1\ra
			= \la\widetilde L_U\Theta u, e_1\ra}.
$$
In particular~$\Im(\Theta)\cap U\ne\{O\}$; otherwise~$\widetilde L_U|_{\Im(\Theta)}$ would be injective, thus~${\Im(\widetilde L_U\Theta)}=\Q^r$, and~$e_1$ could not be orthogonal to~$\Im(\widetilde L_U\Theta)$. Now, it is a well-known fact of linear algebra that given two subspaces~$V, W$ of a vector space~$T$ having the same dimension there exists a subspace~$U$ such that~$T= V\oplus U= W\oplus U$. Hence choosing~$U$ so that~$\Q^n= \Im(N) \oplus U= \Im(\Theta) \oplus U$, we have $\Im(\Theta)\cap U=\{O\}$, and thus the first row of $\widetilde L_U\Theta$ is not identically zero.

Then 
$$
h{1\over m} = (\widetilde L_U\Theta)_1\cdot\pmatrix{\alpha_1 \cr
		   \vdots \cr
		   \alpha_r} + (\widetilde L_U{\bf k})_{1}
$$
and this gives us a contradiction since $\alpha_1,\dots, \alpha_r$ are rationally independent with $1$ by assumption.
%If~$m$ divides~$h$, we have~$[0]$ equals a linear combination with integer coefficients of~$\alpha_1, \dots, \alpha_r$, and since we are assuming~$\alpha_1, \dots, \alpha_r$ to be rationally independent, this can happen only if the first row of~$\widetilde L_U\Theta$ is identically zero, that is the first column of~$\Theta^T\widetilde L_U^T$ is identically zero. This means that the first vector~$e_1$ of the standard basis belongs to~$\ker(\Theta^T\widetilde L_U^T)$, which implies that~$e_1$ is orthogonal to~$\Im(\widetilde L_U\Theta)$, since for any~$u\in\Q^r$, we have
%$${0= \la u, \Theta^T\widetilde L_U^T e_1\ra
%			= \la\Theta u, \widetilde L_U^T e_1\ra
%			= \la\widetilde L_U\Theta u, e_1\ra}.$$
%In particular~$\Im(\Theta)\cap U\ne\{0\}$; otherwise~$\widetilde L_U|_{\Im(\Theta)}$ would be injective, then~${\Im(\widetilde L_U\Theta)}=\Q^r$ and~$e_1$ could not be orthogonal to~$\Im(\widetilde L_U\Theta)$. Since it is a well-known fact of linear algebra that given two subspaces~$V, W$ of a vector space~$T$ having the same dimension there exists a subspace~$U$ such that~$T= V\oplus U= W\oplus U$, we can choose~$U$ such that~$\Q^n= \Im(N) \oplus U= \Im(\Theta) \oplus U$, so~$\Im(\Theta)\cap U=\{0\}$ and this gives us a contradiction. 
\qed

%
%there is a~$r\times r$ submatrix~$N_{r,r}$ of~$N$ with~$\det(N_{r,r})\ne 0$ and let~$\Theta_{r,r}$ the corresponding submatrix in~$\Theta$, i.e.,
%$$N_{r,r}\cdot\pmatrix{\beta_1 \cr
%		   \vdots \cr
%		   \beta_r} = \Theta_{r,r}\cdot\pmatrix{\alpha_1 \cr
%		   \vdots \cr
%		   \alpha_r}.$$
%Then we have
%$$\det(N_{r,r})\pmatrix{\beta_1 \cr
%		   \vdots \cr
%		   \beta_r} = N_{r,r}^*\Theta\cdot\pmatrix{\alpha_1 \cr
%		   \vdots \cr
%		   \alpha_r},$$
%where~$N_{r,r}^*=\det(N_{r,r})\cdot N_{r,r}^{-1}$ and~$N_{r,r}^*\Theta$ is then a~$r\times r$ matrix with integer coefficient, which implies that~$\det(N_{r,r})\beta_1$ is a linear combination with integer coefficients of~$\alpha_1, \dots, \alpha_r$. Since~$\beta_1=1/m$, this means that~$\alpha_1, \dots, \alpha_r$ are rationally dependent contradicting the hypothesis. \qed

\sm We have then two cases to deal with: the rationally independent with $1$ case, and the rationally dependent with $1$ case.

\sm\defin{PureTorsion} Let~$[\phe]\in(\C/\Z)^n$ be of toric degree~$1\le r\le n$. We say that $[\phe]$ is in the {\it torsion-free case}, or simply $[\phe]$ {\it is torsion-free}, if its~$r$-tuples of toric vectors have toric coefficients rationally independent with $1$.

\sm A notion of torsion-free germ of biholomorphism was firstly introduced by \'Ecalle in [\'E]. We shall show in the next section that our notion is equivalent to his; our approach however gives more information on the normalization problem.

\sm We end this section with a couple of results showing how to compute the toric degree, starting with toric degree $1$.

\sm\thm{PrGrado1}{Proposition} {\sl Let $[\phe]\in(\C/\Z)^n$. Then:{\parindent=30pt
\sm\item{(i)} $[\phe]$ has toric degree $1$ with rational toric coefficient if and only if it belongs to $(\Q/\Z)^n$; 
\sm\item{(ii)} $[\phe]$ has toric degree $1$ with toric coefficient in $\C\setminus\Q$ if and only if $[\phe]\not\in(\Q/\Z)^n$, and there exists $\theta\in\Z^n\setminus\{O\}$, with~$\theta_k=0$ if $[\phe_k]=[0]$, such that there is $j_0\in\{1,\dots,n\}$ so that{\parindent=45pt
\sm\item{(a)} $[\phe_{j_0}]\not\in(\Q/\Z)^n$ and
$$
\theta_k[\phe_{j_0}]-\theta_{j_0}[\phe_k]= [0]\tag{EqGr1}
$$
for any $k$ so that $[\phe_k]\ne[0]$; and
\sm\item{(b)} for any representatives $\phe_k$ of $[\phe_k]$, the integer vector~$\phe_{j_0}\theta-\theta_{j_0}\phe$ belongs to the subspace $\hbox{Span}_\Z\{\widehat\theta, -\theta_{j_0} e_1,\dots, \widehat{-\theta_{j_0}e_{j_0}},\dots,-\theta_{j_0} e_{n}\}$, where $\widehat\theta=\theta-\theta_{j_0}e_{j_0}$.\sm} \sm}}

\sm\proof (i) If $\alpha=p/q\in\Q$ then
$$
[\phe]=\left[{p\over q}\theta\right],
$$
hence $[\phe]\in(\Q/\Z)^n$. 

Conversely, if $[\phe_j]= [p_j/q_j]$ with $p_j/q_j\in \Q$ for $j=1, \dots, n$, then, considering $q=q_1\cdots q_n$ we get
$$
[\phe] = \left[\matrix{{p_1q_2\cdots q_n\over q}\cr
				\vdots\cr
				{p_nq_1\cdots q_{n-1}\over q}}\right]
	 =\left[{1\over q} \theta\right], 
$$
and we are done.

\sm(ii) If 
$$
[\phe]= \left[\alpha\pmatrix{\theta_1\cr
							\vdots\cr
							\theta_n}\right],
$$
with $\alpha\in\C\setminus \Q$ and $\theta \in \Z^n\setminus\{O\}$ then it is immediate to verify that $[\phe]\not\in(\Q/\Z)^n$, and~$\theta$ satisfies (a). By assumption, once we choose arbitrarily representatives $\phe_k$ of $[\phe_k]$ we can write~$\phe_k=\alpha\theta_k+m_k$ for suitable $m_k\in\Z$. Then
$$
\theta_k\phe_j-\theta_j\phe_k=\theta_k(\alpha\theta_j+m_j)-\theta_j
(\alpha\theta_k+m_k)=\theta_km_j-\theta_jm_k,
$$
for any $j$ and $k$, thus (b) is verified.

Conversely, let $\theta\in\Z^n\setminus\{O\}$ satisfy the hypotheses. By assumption $[\phe]\not\in(\Q/\Z)^n$ and there is $j_0\in\{1,\dots, n\}$ such that $[\phe_{j_0}]\not\in(\Q/\Z)^n$ satisfies (a) and (b); for the sake of simplicity, we may assume, without loss of generality, $j_0=1$. Let us choose a representative $\phe$ of $[\phe]$ and set 
$$
\theta_j\phe_1-\theta_1\phe_j=k_j\in\Z
$$
for $j=2,\ldots,n$. Condition (b) means that we can find $m_1,\ldots,m_n\in\Z$ so that
$$
\pmatrix{\theta_2 & -\theta_1 & & \cr
    \vdots  & &\ddots &  \cr
    	\theta_n & & & -\theta_1}\cdot
		\pmatrix{m_1\cr
				 \vdots\cr
				 m_n} = \pmatrix{k_2\cr
				 \vdots\cr
				 k_n},\tag{sistemagrado1}
$$
that is
$$
k_j=\theta_j m_1-\theta_1 m_j.
$$
Now we put
$$
\alpha={\phe_1-m_1\over\theta_1}\notin\Q.
$$
Then $[\phe]=[\alpha\theta]$; indeed
$$
\alpha\theta_j={\theta_j(\phe_1-m_1)\over\theta_1}={\theta_j\phe_1-k_j-
\theta_1 m_j\over\theta_1}=\phe_j-m_j.
$$
\qed

\sm\thm{ReGrado1bis}{Remark} Condition (b) of the previous Proposition is necessary. In fact, if we just assume that condition (a) holds, then it is always possible to solve \rf{sistemagrado1} in~$\Q$, but this does not imply that it is solvable in $\Z$. For example the vector
$$
[\phe]=\left[\matrix{{(2i +1)/3}\cr
					i\cr
					{(11 + 10i)/6}}\right]
$$
has toric degree $2$, but if we consider
$$
\theta= \pmatrix{ 2\cr
				3\cr
				5}
$$
we get condition (a) for $j=1$. Moreover, choosing $((2i +4)/3, i, (11 + 10i)/ 6)$ as representative of $[\phe]$, we get 
$$
\pmatrix{k_2\cr
		 k_3} =
\pmatrix{1\cr
		-2}
$$
and it is not difficult to verify that
$$
\pmatrix{3 & -2 & 0 \cr
    	 5 & 0 & -2}\cdot
		\pmatrix{m_1\cr
				 m_2\cr
				 m_3} = \pmatrix{1\cr
				 				 -2}
$$
has no solution $(m_1,m_2, m_3)\in\Z^3$.

%Un esempio lievemente diverso che soddisfa la condizione (a) originale, che era per ogni j,k tali che [\phe_j][\phe_k]\ne[0], senza soddisfare la (b) e' il seguente: \phe=((2i+1)/4, i) e \theta=(2,4).

\sm\thm{Ex1.1}{Example} The vector of~$(\C/\Z)^3$
$$
[\phe_1]= \left[\matrix{{(\sqrt{2} + i)/ 6} \cr
		   {(\sqrt{2} + i)/ 3} \cr
		   {5(\sqrt{2} + i)/ 6} }\right]
$$
has toric degree $1$, since it can be written as
$$
[\phe_1]= \left[{\sqrt{2} + i\over 6}\pmatrix{ 1\cr
		   2 \cr
		   5 }\right].
$$
%Also the vector 
%$$
%[\phe_2] = \left[\matrix{{1/ 6} \cr
%		   {1/3} \cr
%		   {5/6} }\right]
%$$
%has toric degree $1$, since it is equal to
%$$
%[\phe_2]= \left[{1\over 6}\pmatrix{ 1\cr
%		   2 \cr
%		   5 }\right].
%$$

\sm In general, to compute the toric degree of a vector one starts from the trivial representation of Remark \rf{ReToric}, and then uses (the proof of) Lemma \rf{Le1.0} to obtain rationally independent toric coefficients and toric vectors. Then the toric degree is computed as follows (see also Proposition~5.5)
 
\sm\thm{PrCriterio}{Proposition} {\sl Let $\alpha_1, \dots, \alpha_r$ be $1\le r\le n$ rationally independent complex numbers and let $\theta^{(1)}, \dots, \theta^{(r)}\in\Z^n$ be $\Q$-linearly independent integer vectors.
Then:{\parindent=30pt
\sm\item{(i)} if $\alpha_1, \dots, \alpha_r$ are rationally independent with $1$, then $[\phe]=\left[\sum_{k=1}^r \alpha_k \theta^{(k)}\right]$
has toric degree $r$;
\sm\item{(ii)} if $\alpha_1, \dots, \alpha_r$ are rationally dependent with $1$, then $[\phe]=\left[\sum_{k=1}^r \alpha_k \theta^{(k)}\right]$
has toric degree~$r-1$ or $r$.\sm}
}

\sm\proof (i) Let $\alpha_1, \dots, \alpha_r$ be rationally independent with $1$. The toric degree of $[\phe]$ is not greater than $r$. Let us suppose by contradiction that $[\phe]$ has toric degree $s<r$. Then there exist $\eta^{(1)}, \dots, \eta^{(s)}\in\Z^n$ and $\beta_1, \dots, \beta_{s}\in\C$ rationally independent such that
$$
\left[\sum_{k=1}^r\alpha_k\theta^{(k)}\right] = \left[\sum_{k=1}^{s}\beta_k\eta^{(k)}\right]. 
$$
Let~$N$ be the matrix with columns~$\eta^{(1)},\dots, \eta^{(s)}$, and $\Theta$ the matrix with columns~$\theta^{(1)},\dots, \theta^{(r)}$. We have
$$
[\phe]= \left[N\cdot\pmatrix{\beta_1 \cr
		   \vdots \cr
		   \beta_{s}}\right] = \left[\Theta\cdot\pmatrix{\alpha_1 \cr
		   \vdots \cr
		   \alpha_r}\right],
$$
that is there exists an integer vector~${\bf k}\in\Z^n$ such that
$$
N\cdot\pmatrix{\beta_1 \cr
		   \vdots \cr
		   \beta_{s}} = \Theta\cdot\pmatrix{\alpha_1 \cr
		   \vdots \cr
		   \alpha_r} +{\bf k}.
$$		   
Since~$\Theta$ has maximal rank~$r$, the linear map~$\Theta\colon\Q^r\to\Q^n$ is injective and, for every~$U\subseteq\Q^n$ such that~$\Q^n= \Im(\Theta) \oplus U$, there is a linear map~$L_U\colon\Q^n\to\Q^r$ such that~$\ker(L_U) = U$ and~$L_U \Theta= \id$; hence there is a linear map~$\widetilde L_U\colon\Z^n\to\Z^r$ such that~$\widetilde L_U \Theta = h\id$, with~$h\in\Z\setminus\{0\}$.
Then
$$
\widetilde L_U N\cdot\pmatrix{\beta_1 \cr
		   \vdots \cr
		   \beta_{s}}-\widetilde L_U{\bf k}= h\pmatrix{\alpha_1 \cr
		   \vdots \cr
		   \alpha_r}.
$$
Now, $\dim(\ker (\widetilde L_U N)^T) \ge 1$. In particular, there exists $\xi\in\Z^r\setminus \{O\}$ such that $(\widetilde L_U N)^T\xi = O$, that is $\xi^T\widetilde L_U N = O$; therefore
$$
\Z\ni-\xi^T\widetilde L_U{\bf k}
= h\xi^T\pmatrix{\alpha_1 \cr
		   \vdots \cr
		   \alpha_r}		   
= h\la \xi, \alpha\ra,
$$
which is an absurdum, because $\alpha_1, \dots, \alpha_r, 1$ are rationally independent.

\sm (ii) Now we have $\alpha_1, \dots, \alpha_r$ rationally dependent with $1$, and, arguing as in the proof of Lemma \rf{Pr1.0}, we can suppose, without loss of generality, $\alpha_1 =1/m$ and $\alpha_2, \dots, \alpha_r$ rationally independent with $1$. If $m$ divides $\theta^{(1)}_1,\dots, \theta^{(1)}_n$, then $[\phe]= [\sum_{k=2}^r \alpha_k\theta^{(k)}]$ has toric degree $r-1$ thanks to (i). Otherwise, we may assume, without loss of generality, $m, \theta^{(1)}_1,\dots, \theta^{(1)}_n$ coprime.
The toric degree of $[\phe]$ is not greater than~$r$. Let us suppose that $[\phe]$ has toric degree $s<r$. Then there exist $\eta^{(1)}, \dots, \eta^{(s)}\in\Z^n$ and $\beta_1, \dots, \beta_{s}\in\C$ such that
$$
\left[\sum_{k=1}^r\alpha_k\theta^{(k)}\right] = \left[\sum_{k=1}^{s}\beta_k\eta^{(k)}\right], 
$$
%If $\beta_1,\dots, \beta_s$ are rationally independent with $1$, 
thus we have
$$
[m\phe] =\left[\sum_{k=2}^r \alpha_k \cdot m\theta^{(k)} \right]= \left[\sum_{k=1}^s \beta_k \cdot m\eta^{(k)} \right],
$$
and, since $\alpha_2, \dots, \alpha_r$ are rationally independent with $1$, by (i) we get $s=r-1$. 
%If $\beta_1,\dots,\beta_s$ are rationally dependent with $1$, hence we can assume that we have a reduced $s$-tuple and $\beta_1= 1/p$, i.e.,
%$$
%[\phe]=\left[{1\over m} \theta^{(1)} + \sum_{k=2}^r \alpha_k \theta^{(k)} \right]= \left[{1\over p} \eta^{(1)} + \sum_{k=2}^s \beta_k \eta^{(k)} \right].
%$$
%Then
%$$
%[m\phe] =\left[\sum_{k=2}^r \alpha_k \cdot m\theta^{(k)} \right]= \left[{m\over p} \eta^{(1)} + \sum_{k=2}^s \beta_k \cdot m\eta^{(k)} \right],
%$$
%and we have two cases. If $p$ divides $m$ then 
%$$
%\left[{m\over p} \eta^{(1)} + \sum_{k=2}^s \beta_k \cdot m\eta^{(k)} \right]= \left[\sum_{k=2}^s \beta_k \cdot m\eta^{(k)} \right]
%$$
%and we are in the case with $\alpha_2, \dots, \alpha_r$ rationally independent with $1$, therefore $r=s$. If $p$ does not divide $m$, then
%$$
%[mp\phe] =\left[\sum_{k=2}^r \alpha_k \cdot mp\theta^{(k)} \right]= \left[\sum_{k=2}^s \beta_k \cdot mp\eta^{(k)} \right],
%$$ 
%and again, we must have $r=s$.
\qed

\sm\thm{ReCriterio}{Remark} Note that both cases in (ii) can occur. In fact, 
it is not difficult to verify that
$$
[\phe_1]= \left[\matrix{1/2\cr
						\sqrt{2}\cr	
						i}\right]= \left[ {1\over 2} e_1 + \sqrt{2} e_2 + i e_3\right]
$$
has toric degree $3$. However, if we consider,
$$
[\phe_2]= \left[{1\over 2} \pmatrix{1\cr		
									1\cr
									1} + 
									{\sqrt{2}-1\over 2} \pmatrix{1\cr											           1\cr							                   3} +
									{i\over 2} \pmatrix{0\cr											       1\cr							               1}\right],
$$
then
$$
[\phe_2]=\left[\matrix{\sqrt{2}/2\cr
					(\sqrt{2} + i)/2\cr	
						(-2 + 3\sqrt{2} + i)/2}\right]= \left[{\sqrt{2}\over 2}\pmatrix{1\cr
		  1\cr
		  3} + {i\over 2} \pmatrix{0\cr											 1\cr							         1}\right],
$$
so the toric degree is $2$. Proposition~5.5 will show how to distinguish between the two cases of Proposition \rf{PrCriterio}.(ii).

%\xi^T\widetilde L_U N\cdot\pmatrix{\beta_1 \cr
%		   \vdots \cr
%		   \beta_{s}}
%%%%%%%%%%%%%%%%%%%%%%%%%%%%%%%%%%%%%%%%%%%%%%%%%%%%%%%%%%%%%%%%%%%%%%%%%%%%%%%%
\sect Torsion
%%%%%%%%%%%%%%%%%%%%%%%%%%%%%%%%%%%%%%%%%%%%%%%%%%%%%%%%%%%%%%%%%%%%%%%%%%%%%%%%

\sm In [\'E], \'Ecalle introduced the following notion. 

\sm\defin{detorsion} Let $\lambda\in(\C^*)^n$. The {\it torsion} of~$\lambda$ is the natural integer~$\tau$ such that
$$
{1\over \tau}2\pi i\Z = (2\pi i \Q)\cap \left( (2\pi i \Z)\bigoplus_{1\le j\le n}((\log \lambda_j)\Z)\right).\tag{eqecalle}
$$

\sm Translated in our notation, \rf{eqecalle} becomes
$$
{1\over \tau}\Z = \Q\cap \left( \Z\bigoplus_{1\le j\le n}\phe_j\Z\right),
$$
where $\phe$ is a representative of the unique $[\phe]\in(\C/\Z)^n$ such that $\lambda=\exp(2\pi i[\phe])$.

\sm The torsion is well-defined, as the following result shows (and whose proof
describes how to explicitly compute the torsion).

\sm\thm{Le1.4}{Proposition} {\sl The torsion of a~$n$-tuple~$(\lambda_1,\dots, \lambda_n)\in(\C^*)^n$ is a well-defined natural integer. Furthermore, writing $\lambda=e^{2\pi i[\phe]}$, if $[\phe]$ is torsion-free, then $\tau=1$; otherwise $\tau$ divides the denominator of the first toric coefficient in a reduced representation of $[\phe]$.}

\sm\proof Let~$[\phe]\in(\C/\Z)^n$ be the unique vector such that $\lambda=\exp(2\pi i[\phe])$, let~$1\le r\le n$ be its toric degree and let~$\theta^{(1)},\dots, \theta^{(r)}$ be a~$r$-tuple of toric vectors associated to~$[\phe]$ with coefficients~$\alpha_1,\dots, \alpha_r$. 

Our aim is to determine the structure of the set
$$
R=\Q\cap\left(\Z\bigoplus_{j=1}^n \phe_j\Z\right),
$$
that is of the set of rational numbers $x$ that can be expressed in the form
$$
\Q\ni x=m_0+m_1\phe_1+\cdots+m_n\phe_n
$$
with $m_0,\ldots,m_n\in\Z$. Write, as usual,
$$
\phe_j=h_j+\sum_{k=1}^r \alpha_k\theta^{(k)}_j,
$$
with $h_j\in\Z$. Then
$$
\eqalign{
x&=(m_0+m_1h_1+\cdots+m_nh_n)+m_1\sum_{k=1}^r \alpha_k\theta^{(k)}_1+\cdots+
m_n\sum_{k=1}^r \alpha_k\theta^{(k)}_n\cr
&=\widetilde m+\sum_{k=1}^r\alpha_k\langle M,\theta^{(k)}\rangle,
\cr}
$$
where $\widetilde m\in\Z$ and $M\in\Z^n$ are generic. If $\alpha_1,\ldots,\alpha_r$ are rationally independent with~1, it follows that $x\in\Q$ if
and only if $\langle M,\theta^{(1)}\rangle=\cdots=\langle M,\theta^{(r)}
\rangle=0$, and thus $R=\Z$ and $\tau=1$.

If $\alpha_1,\ldots,\alpha_r$ are not rationally independent with~1, let us use instead the reduced representation, with $\beta_1=1/m$, the
remaining coefficients $\beta_2,\ldots,\beta_r$ rationally independent with~1, and with $\eta^{(1)},\ldots,\eta^{(r)}$ as toric vectors. We get
$$
x=\widetilde m+{1\over m}\langle M,\eta^{(1)}\rangle+\sum_{k=2}^r\beta_k
\langle M,\eta^{(k)}\rangle.
$$
Therefore $x\in\Q$ if and only if $\langle M,\eta^{(2)}\rangle=\cdots=\langle M,\eta^{(r)}\rangle=0$, and moreover in that case
$$
x=\widetilde m+{1\over m}\langle M,\eta^{(1)}\rangle.
$$
Now, the set
$$
S=\{\langle M,\eta^{(1)}\rangle\mid M\in\Z^n, \langle M,\eta^{(2)}\rangle
=\cdots=\langle M,\eta^{(r)}\rangle=0\}
$$
is an ideal of~$\Z$; therefore $S=q\Z$ for some $q\in\N$. It follows that
$$
R=\Z\oplus{q\over m}\Z=\Z\oplus{\widetilde q\over\widetilde m}\Z={1\over\widetilde m}\Z,
$$
where $\widetilde q$ and $\widetilde m$ are coprime, and $q/m=\widetilde q/\widetilde m$.
Hence $\tau=\widetilde m$, and we are done. \qed

\sm\thm{Re5.1}{Remark} Note that, in the previous proof, $S\ne\{O\}$, i.e., $q\ne 0$. Indeed, $S=\{O\}$ if and only if the kernel in~$\Z^n$ of the linear form $(\eta^{(1)})^T$ contains the intersection of the kernels in~$\Z^n$ of the linear forms $(\eta^{(2)})^T,\ldots,(\eta^{(r)})^T$. It is easy to see that this implies that the kernel in~$\Q^n$ of the linear form~$(\eta^{(1)})^T$
contains the intersection of the kernels in $\Q^n$ of the linear forms~$(\eta^{(2)})^T,\ldots,(\eta^{(r)})^T$. But this implies that the
linear form $(\eta^{(1)})^T$ is a $\Q$-linear combination of~$(\eta^{(2)})^T,\ldots,(\eta^{(r)})^T$, and so $\eta^{(1)},\ldots,
\eta^{(r)}$ are $\Q$-linearly dependent, impossible.

%% Remark che va nella tesi
%\sm\thm{RemTesi}{Remark} Note that, if 
%$$
%[\phe] = \left[{1\over \tau q} \eta^{(1)} + \sum_{k=2}^r \beta_k\eta^{(k)}\right],
%$$
%with $\tau$ and $q$ coprime, then there exist $\gamma_2,\dots, \gamma_r\in\C^*$ and $\xi^{(1)}, \dots, \xi^{(r)}\in\Z^n$ such that 
%$$
%[\phe] = \left[{1\over \tau} \xi^{(1)} + \sum_{k=2}^r \gamma_k\xi^{(k)}\right].
%$$
%In fact, there exist $a,b\in\Z^*$ such that $a\tau + b q = 1$, then
%$$
%{1\over \tau q} = {a\over q} + {b\over \tau}.
%$$
%Moreover
%$$
%\left[{b\over q} \eta^{(1)} + \sum_{k=2}^r \beta_k\eta^{(k)}\right]
%$$ 
%has no torsion, then there exist $\gamma_2,\dots, \gamma_r\in\C^*$ and $\xi^{(2)}, \dots, \xi^{(r)}\in\Z^n$ such that 
%$$
%\left[{b\over q} \eta^{(1)} + \sum_{k=2}^r \beta_k\eta^{(k)}\right] = \left[\sum_{k=2}^r \gamma_k\xi^{(k)}\right].
%$$
%Hence
%$$
%\eqalign{[\phe] 
%&=\left[{1\over \tau q} \eta^{(1)} + \sum_{k=2}^r \beta_k\eta^{(k)}\right]\cr
%&=\left[\left({a\over q} + {b\over \tau} \right)\eta^{(1)} + \sum_{k=2}^r \beta_k\eta^{(k)}\right]\cr
%&=\left[{b\over \tau} \eta^{(1)} + \sum_{k=2}^r \gamma_k\xi^{(k)}\right],}
%$$
%and we are done. 

\sm The next result explains the terminology of Definition \rf{PureTorsion}.

%shows that we cannot have $\la P, \eta^{(1)}\ra \in m\Z$ for any $P\in\Z^n$ such that $\la P,\eta^{(k)}\ra =0$ for $k=2, \dots, r$.

\sm\thm{PrNoImpureTorsionFree}{Theorem} {\sl Let $\lambda= e^{2\pi i [\phe]}\in(\C^*)^n$. Then~$[\phe]$ is torsion-free if and only if the torsion of $\lambda$ is $1$.}

\sm\proof If $[\phe]$ is torsion-free, then the toric coefficients of a toric $r$-tuple associated to $[\phe]$ are rationally independent with $1$, and the torsion $\tau$ is~$1$, by Proposition \rf{Le1.4}.

Conversely, let~$\eta^{(1)},\dots, \eta^{(r)}$ be a reduced~$r$-tuple of toric vectors associated to~$[\phe]$ with toric coefficients~$1/m,\beta_2,\dots,\beta_r$. Let us assume by contradiction that the torsion $\tau$ of $[\phe]$ is $1$. From the proof of Proposition \rf{Le1.4} it is clear that we have $\tau=1$ if and only if $\la P, \eta^{(1)}\ra \in m\Z$, for any~$P\in\Z^n$ such that $\la P,\eta^{(k)}\ra =0$ for $k=2, \dots, r$.

Since $\eta^{(1)}, \dots, \eta^{(r)}$ are a toric $r$-tuple, we may assume, without loss of generality, that the matrix $A$ of $M_{n\times n}(\Z)$ with columns $\eta^{(2)}, \dots, \eta^{(r)}, e_r, \dots, e_n$ is invertible in $M_{n\times n}(\Q)$. Denote by $N'$ the matrix in $M_{(r-1)\times(r-1)}(\Z)$ 
$$
N'= \pmatrix{\eta^{(2)}_1 & \dots & \eta^{(r)}_1\cr
					\vdots &      &\vdots\cr
					\eta^{(2)}_{r-1} & \dots & \eta^{(r)}_{r-1}},
$$
and by $N''$ the matrix in $M_{(n-r+1)\times(r-1)}(\Z)$
$$
N''= \pmatrix{\eta^{(2)}_{r} & \dots & \eta^{(r)}_{r}\cr
					\vdots &      &\vdots\cr
					\eta^{(2)}_{n} & \dots & \eta^{(r)}_{n}}.
$$
Then 
$$
A=\pmatrix{ N' & O \cr
			N'' & I_{n-r+1}}
$$
and $\det(A)=\det(N')\ne 0$. 

We claim that, up to pass to another toric $r$-tuple $\widehat\eta^{(1)}, \eta^{(2)}, \dots, \eta^{(r)}$, we may assume that~$m=\det(N')$ and $\widehat\eta^{(1)}\in \{0\}^{r-1}\times \Z^{n-r+1}$. In fact, $\eta^{(k)} = A^{-1}e_{k-1}$ for $k=2, \dots, r$, with~$A^{-1}\in M_{n\times n}(\Q)$. Hence~$P\in\Z^n$ is such that $\la P,\eta^{(k)}\ra =0$ for $k=2, \dots, r$ if and only if~$\la A^T P, e_j \ra = 0$ for $j=1, \dots, r-1$, that is $A^T P\in \{0\}^{r-1}\times \Z^{n-r+1}$. Now, we have
$$
A^T P = \pmatrix{ N'^T & N''^T \cr
			O & I_{n-r+1}} \pmatrix{ P'\cr
									 P''} \in \{0\}^{r-1}\times \Z^{n-r+1}
$$
if and only if
$$
P = \pmatrix{-(N'^T)^{-1}N''^T P''\cr
				P''} \quad\hbox{with}\quad P''\in\Z^{n-r+1}~\hbox{and}~(N'^T)^{-1}N''^T P''\in\Z^{r-1},
$$
that is
$$
P''\in\Z^{n-r+1}~\hbox{and}~(N'^+)^{T}N''^T P''\in \det(N')\Z^{r-1}
$$
where $(N'^+)^{T}\in M_{(r-1)\times (r-1)}(\Z)$ and $(N'^+)^{T}N'=\det(N')I_{r-1}$. In particular, since we are assuming
$$
\la P,\eta^{(k)}\ra= 0~\hbox{for}~k=2,\dots, r \,\Longrightarrow\, \la P,\eta^{(1)}\ra\in m \Z,\tag{ipotesiassurdo}
$$
in particular, we get 
$$
\la A^T P,A^{-1}\eta^{(1)}\ra= \left\la\pmatrix{O\cr P''} ,A^{-1}\eta^{(1)}\right\ra\in m\Z
$$
for any $P''\in \det(N')\Z^{n-r+1}$. Then there exist $q_1, \dots, q_{r-1}\in \Q$ and $\widehat \eta^{(1)}\in \{0\}^{r-1}\times \Z^{n-r+1}$ such that 
$$
A^{-1} \eta^{(1)} = q_1 e_1+ \cdots + q_{r-1}e_{r-1} + {m\over\det(N')} \widehat \eta^{(1)},
$$
that is
$$
\eta^{(1)} = q_1 \eta^{(2)}+ \cdots + q_{r-1}\eta^{(r)} + {m\over\det(N')} \widehat \eta^{(1)},
$$
thus we get
$$
\eqalign{
[\phe]&=\left[{1\over m} \eta^{(1)} + \sum_{k=2}^r \beta_k\eta^{(k)}\right]\cr
	 &=\left[{1\over m} {m\over\det(N')}\widehat\eta^{(1)} + \sum_{k=2}^r \left(\beta_k+ {q_{k-1}\over m}\right)\eta^{(k)}\right]\cr
	 &=\left[{1\over \det(N')} \widehat\eta^{(1)} + \sum_{k=2}^r \widetilde\beta_k\eta^{(k)}\right].}
$$
Note that $\widetilde \beta_2,\dots, \widetilde \beta_r$ are rationally independent with $1$.

Now we can assume that \rf{ipotesiassurdo} holds with $m=\det(N')$ and $\eta^{(1)}\in \{0\}^{r-1}\times \Z^{n-r+1}$. We claim that there exist $\gamma_2, \dots,\gamma_r\in\C^*$ such that $[\phe]=[\sum_{k=2}^r \gamma_k\eta^{(k)}]$, i.e., $[\phe]$ has toric degree~$r-1$, contradicting the hypotheses. We can have $[\phe]=[\sum_{k=2}^r \gamma_k\eta^{(k)}]$ with $\gamma_2, \dots,\gamma_r\in\C^*$, if there exists $\theta'\in \Z^{r-1}$ such that
$$
\pmatrix{\gamma_2\cr
		\vdots\cr
		\gamma_r} = \pmatrix{\beta_2\cr
		\vdots\cr
		\beta_r} + N'^{-1}\theta',
$$
and $\theta'\in\Z^{r-1}$ is a solution
$$
N''N'^+\pmatrix{x_1\cr
				\vdots\cr
				x_{r-1}} \equiv \pmatrix{\eta^{(1)}_r\cr
										\vdots\cr
										\eta^{(1)}_n}~~\hbox{mod}~m\Z^{n-r+1}.\tag{sistemaassurdo}
$$
In fact, since $N'' N'^{-1} = (1/m) N'' N'^+$, this implies 
$$
{1\over m}\eta^{(1)} = {1\over m}\pmatrix{O\cr
								\eta''^{(1)}}\equiv \pmatrix{O\cr
								N'' N'^{-1}\theta'},\tag{eqstella}
$$		
modulo $\Z$, where $\eta''^{(1)}=(\eta^{(1)}_r,\dots,\eta^{(1)}_n)$, hence
$$
\eqalign{
[\phe] &= \left[{1\over m} \eta^{(1)} + N\pmatrix{\beta_2\cr
		                           \vdots\cr
		                           \beta_r}\right]\cr
	  &= \left[\pmatrix{O\cr
					  N'' N'^{-1}\theta'}+ \pmatrix{N'\cr
								        N''}\left(\pmatrix{\gamma_2\cr
		\vdots\cr
		\gamma_r}- N'^{-1}\theta'\right)\right]\cr
	  &= \left[\pmatrix{O\cr
					  N'' N'^{-1}\theta'}+ \pmatrix{N'\cr
								        N''}\pmatrix{\gamma_2\cr
								               \vdots\cr
								               \gamma_r}-\pmatrix{\theta'\cr	
     N''N'^{-1}\theta'}\right]\cr
    &=\left[N\pmatrix{\gamma_2\cr
	                \vdots\cr	
	                \gamma_r}\right].}
$$

Now we prove that, if \rf{ipotesiassurdo} holds with $m=\det(N')$ and $\eta^{(1)}\in \{0\}^{r-1}\times \Z^{n-r+1}$, then there exists a solution $\theta'\in \Z^{r-1}$ of \rf{sistemaassurdo}. In fact, if $P''\not\in m\Z^{n-r+1}$ is a multi-index such that~$P''^T N'' N'^+ \in m\Z^{r-1}$, then by \rf{ipotesiassurdo} we have $P''^T \eta''^{(1)}\in m\Z$, where we use the same notation of \rf{eqstella}; thus, since up to reorder the indices we may assume that the last coordinate of $P''$ is not in~$m\Z$, we can substitute $P''^T N'' N'^+ x \equiv P''^T \eta''^{(1)}$ to the last equation of \rf{sistemaassurdo}, and we have to solve a system with one equation less. We iterate this procedure for a set of generators of a complement of $m\Z^{n-r+1}$ in the lattice of $P''$ until, up to reordering, we get
$$
B \pmatrix{x_1\cr
			\vdots\cr
		  x_{r-1}} \equiv \pmatrix{\eta^{(1)}_r\cr
										\vdots\cr
										\eta^{(1)}_{r+h-1}}~~\hbox{mod}~m\Z^{h}
$$
where $1\le h\le n-r+1$, $B\in M_{h\times (r-1)}(\Z)$ is the matrix of the first $h$ rows of $N'' N''^+$, and for any $R\not\in m\Z^h$, we have $R^T B\not \in m\Z^{r-1}$, that is $B$ has maximal rank modulo $m$.

If $h=1$, then we have
$$
b_1 x_1 + \cdots + b_{r-1} x_{r-1} \equiv \eta^{(1)}_1~~\hbox{mod}~m\Z.\tag{eqassurdo}
$$ 
If $b_1, \dots, b_{r-1}, m$ are coprime it is obvious that \rf{eqassurdo} is solvable. If the greatest common divisor of $b_1, \dots, b_{r-1}, m$ is $p>1$, then $m= q p$ and $q (b_1, \dots, b_{r-1})\in m\Z^{r-1}$, hence, by \rf{ipotesiassurdo}, we must have $\eta^{(1)}_1\in p\Z$ too, thus 
$$
{b_1\over p} x_1 + \cdots + {b_{r-1}\over p} x_{r-1} \equiv {\eta^{(1)}_1\over p}~~\hbox{mod}~{m\over p}\Z
$$
is solvable.

Let us now suppose $1< h\le n-r+1$. Since $B$ has maximal rank modulo $m$, there exists~$B^+$ in~$M_{(r-1)\times h}(\Z)$ such that $B^+ B \equiv d I_{r-1}$, modulo $m\Z$ where $d\ne m$. 
Thus we have
$$
d\pmatrix{x_1\cr
			\vdots\cr
		  x_{r-1}} \equiv B^+ \pmatrix{\eta^{(1)}_r\cr
										\vdots\cr
										\eta^{(1)}_{r+h-1}}~~\hbox{mod}~m\Z^{h}.
$$
If $d$ and $m$ are coprime, we are done. Otherwise, let $p$ be greatest common divisor of $d$ and~$m$, and let $q=m/p$. Since $B^+ B \equiv d I_{r-1}$ modulo $m\Z$, we have $q B^+ B \equiv O$ modulo $m\Z$, thus, since we are assuming that for any $R\not\in m\Z^h$, we have $R^T B\not \in m\Z^{r-1}$, it has to be $q B^+\equiv O$ modulo $m\Z$, that is $B^+\equiv p\widetilde B$ modulo $m\Z$. 
Therefore we have
$$
{d\over p}\pmatrix{x_1\cr
			\vdots\cr
		  x_{r-1}} \equiv \widetilde B \pmatrix{\eta^{(1)}_r\cr
										\vdots\cr
										\eta^{(1)}_{r+h-1}}~~\hbox{mod}~{m\over p}\Z^{h},
$$
which is solvable, as we wanted.\qed

\me The torsion case is more delicate and difficult to deal. First, given $[\phe]\in(\C/\Z)^n$ with toric degree $1\le r\le n$ and torsion $\tau\ge2$, and a reduced toric $r$-tuple $\eta^{(1)},\dots, \eta^{(r)}$, we have
$$ 
\bigcap_{k=2}^r \res_j^+(\eta^{(k)}) \supseteq \res_j([\phe])\supseteq
\bigcap_{k=1}^r \res_j^+(\eta^{(k)}),\tag{trecasiintro},
$$
yielding a subdivision in more subcases, all realizable (we have examples for all of them) and, surprisingly, having very different behaviours ones from the others; we have cases similar to the case of germs of vector fields (even if we have torsion!), and cases that are indeed different.
In particular, considering iterates of $f$ to reduce to the torsion-free case hides very interesting phenomena, and it does not allow to see that some torsion cases can be directly studied. 
Moreover, we have explicit (and computable) techniques to decide in which subcase a given $[\phe]\in(\C/\Z)^n$ is.

\sm\thm{ExPurePresenceTorsion1}{Example} Let us consider the vector
$$
[\phe]= \left[{1\over 6}\pmatrix{1 \cr
		   3 } + \sqrt{2}\pmatrix{1 \cr
		   -6}
		   \right]\in(\C/\Z)^2,
$$
of toric degree $2$.
%In this case $\ca D(1/6,\sqrt{2}) = \{(6h, 0)\mid h\in\Z\}$. 
We have
$$
\la P,\eta^{(2)}\ra = p_1 - 6 p_2 = 0
$$
if and only if
$$
P\in \pmatrix{ 6\cr
				1}\Z,
$$
hence
$$
\res_1^+(\eta^{(2)}) = \{(6h+1,h)\mid h\in\N\setminus\{0\}\} \quad\hbox{and}\quad \res_2^+(\eta^{(2)}) = \{(6h, h+1)\mid h\in\N\setminus\{0\}\},
$$
and
$$
\la (6h,h), \eta^{(1)}\ra \in 9\Z,
$$
that is
$$
S=9\Z,
$$
%hence
%$$
%\ca D (1/6, \sqrt{2}) \cap {\rm Adm}_j(\eta^{(1)}, \eta^{(2)})\ne \{O\}
%$$
%for $j=1,2$ 
and the torsion is clearly $2$. Moreover, we have
$$
[\phe]= \left[{1\over 2}\pmatrix{1 \cr
		   1 } + {3\sqrt{2}-1\over 3}\pmatrix{1 \cr
		   -6}
		   \right]\in(\C/\Z)^2.
$$

\sm Using the torsion $\tau$ of a vector, we obtain a complete criterion to compute the toric degree of a vector, as next result shows.
 
\sm\thm{PrCriterio2}{Proposition} {\sl Let $[\phe]\in(\C/\Z)^n$ and let $\tau$ be its torsion. If 
$$
[\phe]= \left[{1\over m} \eta^{(1)} + \sum_{k=2}^r\beta_k\eta^{(k)}\right],
$$
with $\eta^{(1)}\not \in m\Z^n$, then $[\phe]$ has toric degree $r$ if and only
if the torsion of $[\phe]$ is $\tau>1$, the coefficients $\beta_2,\dots, \beta_r$ are rationally independent with $1$, and the integer vectors $\eta^{(1)}, \dots, \eta^{(r)}$ are $\Q$-linearly independent.}

\sm\proof It follows from Lemma \rf{Pr1.0}, Proposition \rf{Le1.4} and from the proof of Theorem \rf{PrNoImpureTorsionFree}. 

%%%%%%%%%%%%%%%%%%%%%%%%%%%%%%%%%%%%%%%%%%%%%%%%%%%%%%%%%%%%%%%%%%%%%%%%%%%%%%%%
\sect Poincar\'e-Dulac Normal Form in the torsion-free case
%%%%%%%%%%%%%%%%%%%%%%%%%%%%%%%%%%%%%%%%%%%%%%%%%%%%%%%%%%%%%%%%%%%%%%%%%%%%%%%%

%\sm{\sf to be checked}\thm{Re1.3}{Remark} {\sl If~$\phe\in\C^n$ has toric degree~$r>1$ then the~$r$-tuple of toric vectors associated to~$\phe$ is not necessarily unique.} Let us consider, for example
%$$\phe = (3+4 i, 2+6 i, -1 + 2 i).$$ 
%The toric degree cannot be~$1$, since it is immediate to verify that~$\phe$ cannot be written as the product of a complex number times an integer vector.
%The toric degree is in fact~$2$, since choosing~$\alpha_1=i, \alpha_2= 1+i, \theta^{(1)}= (1, 4, 3)$ and~$\theta^{(2)}= (3, 2, -1)$, we have
%$$\eqalign{\alpha_1 \theta^{(1)} + \alpha_2 \theta^{(2)}&=i \left(1, 4, 3\right) + (1+i) \left(3, 2, -1\right) \cr&= (i+3 +3i,4i +2+2i, 3i -1-i)\cr&= \phe}$$ 
%and~$\alpha_1,\alpha_2$ are rationally independent. Anyway, choosing~$\beta_1=(-3 + 16 i)/6$, ${\beta_2= (-3-4i)/6}$, $\gamma^{(1)}= (0,1,1)$ and~$\gamma^{(2)} = (-6, -5 ,1)$, it is also true
%$$\eqalign{\phe&=\beta_1 \gamma^{(1)}+ \beta_2 \gamma^{(2)}\cr&= {-3 + 16 i\over 6} \left(0, 1, 1\right) + {-3-4i\over 6} \left(-6, -5, 1\right) \cr &= \left(-6{(-3-4i)\over 6}, {-3 + 16 i+ 15+ 20i\over 6}, {-3 + 16 i-3-4i\over 6} \right)}$$
%and~$\beta_1,\beta_2$ are rationally independent. 

\sm In the torsion-free case, it is not difficult to show that we can compute the resonances of~$[\phe]$, which are multiplicative, using the additive resonances of one of its associated~$r$-tuples of toric vectors, as the next result shows.

\sm\thm{Le1.2}{Lemma} {\sl Let~$[\phe]\in(\C/\Z)^n$ be of toric degree~$1\le r\le n$ and torsion-free. Then for any~$r$-tuple of toric vectors,~$\theta^{(1)},\dots, \theta^{(r)}$, associated to~$[\phe]$ we have
$$
\res_j([\phe]) = \bigcap_{k=1}^r \res_j^+(\theta^{(k)})
$$
for every~$j=1, \dots, n$.}

\sm\proof We have
$$
[\la Q,\phe \ra - \phe_j] = \left[\sum_{k=1}^r \alpha_k \left(\la Q,\theta^{(k)} \ra - \theta^{(k)}_j\right)\right] \tag{eq1.2}
$$
and, since~$\alpha_1,\dots, \alpha_r$ are rationally independent with $1$, the right-hand side of \rf{eq1.2} vanishes if and only if~$\la Q,\theta^{(k)} \ra - \theta^{(k)}_j=0$ for every~$k=1,\dots, r$. \qed

\sm\thm{ExPureTorsion1}{Example} Let us consider the torsion-free vector
$$
[\phe]= \left[\sqrt{2}\pmatrix{3 \cr
		   2 \cr
		   -1 } + 2i\pmatrix{2 \cr
		   3 \cr
		   1 }
		   \right]\in(\C/\Z)^3,
$$
of toric degree $2$.
Then
$$
\cases{\displaystyle \la P,\theta^{(1)}\ra = 3 p_1 + 2 p_2 - p_3 = 0 \cr\noalign{\sm}
\displaystyle \la P,\theta^{(2)}\ra = 2 p_1 + 3 p_2 + p_3 = 0}
$$
for some $P\in\Z^n$, if and only if
$$
P\in\pmatrix{1\cr
			-1\cr 
			1}\Z.
$$
Hence in this case
$$
\res_1([\phe]) =\res_3([\phe]) = \void~~\hbox{and}~~\res_2([\phe]) = \{(1,0,1)\}.
$$

\sm\thm{ExPureTorsion2}{Example} Let us consider the vector
$$
[\phe]= \left[\sqrt{2}\pmatrix{3 \cr
		   2 \cr
		   -1 } + 2i\pmatrix{2 \cr
		   -3 \cr
		   1 }
		   \right]\in(\C/\Z)^3.
$$
Again, $[\phe]$ has toric degree $2$ and it is torsion-free.
In this case, we have
$$
\cases{\displaystyle \la P,\theta^{(1)}\ra = 3 p_1 + 2 p_2 - p_3 = 0 \cr\noalign{\sm}
\displaystyle \la P,\theta^{(2)}\ra = 2 p_1 - 3 p_2 + p_3 = 0}
$$
for some $P\in\Z^n$, if and only if
$$
P\in\pmatrix{1\cr
			5\cr 
			13}\Z.
$$
Hence
$$
\eqalign{
&\res_1([\phe]) =\left\{(q+1, 5 q, 13 q)\mid q\in\N\setminus\{0\}\right\}\cr
&\res_2([\phe]) = \left\{(q, 5 q +1, 13 q)\mid q\in\N\setminus\{0\}\right\}\cr
&\res_3([\phe]) = \left\{(q, 5 q, 13 q +1 )\mid q\in\N\setminus\{0\}\right\}.}
$$

\sm We have the following immediate corollary of Lemma \rf{Le1.2}.

\sm\thm{Co1.1}{Corollary} {\sl Let~$\lambda\in(\C^*)^n$ and let~$[\phe]\in(\C/\Z)^n$ be such that~$\lambda= e^{2\pi i [\phe]}$. If~$[\phe]$ is torsion-free and has toric degree~$1\le r\le n$, then for every~$r$-tuple~$\theta^{(1)}, \dots, \theta^{(r)}$ of toric vectors associated to~$[\phe]$ we have
$$
{\rm Res}_j(\lambda) =\bigcap_{k=1}^r \res_j^+(\theta^{(k)})
$$
for every~$j=1, \dots, n$.}

\sm\thm{Le1.1}{Lemma} {\sl Let~$[\phe]\in(\C/\Z)^n$ be of toric degree~$1\le r\le n$ and torsion-free. Then for any~$r$-tuple of toric vectors,~$\theta^{(1)},\dots, \theta^{(r)}$, associated to~$[\phe]$ we have $\theta^{(k)}_j = \theta^{(k)}_h$ whenever $[\phe_j]=[\phe_h]$, for every $k=1,\dots, r$.}

\sm\proof If~$[\phe_j]=[\phe_h]$, then
$$
\left[\alpha_1\theta^{(1)}_j + \cdots + \alpha_r\theta^{(r)}_j\right] =\left[\alpha_1\theta^{(1)}_h + \cdots + \alpha_r\theta^{(r)}_h\right];
$$
hence there exists $m\in\Z$, such that
$$
\alpha_1\left(\theta^{(1)}_j-\theta^{(1)}_h \right) + \cdots + \alpha_r\left(\theta^{(r)}_j-\theta^{(r)}_h\right)=m,
$$
and, since~$\theta^{(k)}_j-\theta^{(k)}_h\in\Z$ for~$k=1,\dots, r$, the assertion follows from the rational independence with $1$ of~$\alpha_1,\dots, \alpha_r$. \qed 

%\sm\defin{DeTorsionFree} Let~$f$ be germ of biholomorphism of~$\C^n$ fixing the origin~$O$ with~$\{\lambda_1, \dots, \lambda_n\}$ as spectrum of~$\d f_O$. We say that $f$ is {\it torsion-free} if the unique~$[\phe]=([\phe_1],\dots, [\phe_n])\in(\C/\Z)^n$ such that~$\lambda=e^{2\pi i [\phe]}$ is torsion-free. 

\sm\defin{PureGerm} Let~$f$ be a germ of biholomorphism of~$\C^n$ fixing the origin. We say that $f$ is {\it torsion-free} if, denoted by $\lambda=\{\lambda_1, \dots, \lambda_n\}$ the spectrum of~$\d f_O$, the unique $[\phe]\in(\C/\Z)^n$ such that~$\lambda=e^{2\pi i [\phe]}$ is in the torsion-free case.

\sm We have then the following complete description of Poincar\'e-Dulac holomorphic normalization in the torsion-free case.

\sm\thm{Te1.1}{Theorem} {\sl Let~$f$ be a germ of biholomorphism of~$\C^n$ fixing the origin~$O$, of toric degree~$1\le r\le n$ and in the torsion-free case. Then~$f$ admits a holomorphic Poincar\'e-Dulac normalization if and only if there exists a holomorphic effective action on~$(\C^n,O)$ of a torus of dimension~$r$ commuting with~$f$ and such that the columns of the weight matrix of the action are a~$r$-tuple of toric vectors associated to~$f$.}

\sm\proof It follows from Theorem \rf{Te2.1}, Lemma \rf{Le1.1} and Corollary \rf{Co1.1}. \qed

%Let us suppose~$f$ is in Poincar\'e-Dulac holomorphic normal form and let~$\theta^{(1)},\dots, \theta^{(r)}$ be a~$r$-tuple of toric vectors associated to~$f$. Then it commutes with the linear effective action of~$\T^r$ defined by
%$$A(x, z) = \diag\left(e^{2\pi i \sum_{k=1}^r x_k \theta^{(k)}_j}\right)z.$$ 
%The converse follows from Theorem \rf{Te2.1} and Corollary \rf{Co1.1}. \qed

%%%%%%%%%%%%%%%%%%%%%%%%%%%%%%%%%%%%%%%%%%%%%%%%%%%%%%%%%%%%%%%%%%%%%%%%%%%%%%%%
\sect Poincar\'e-Dulac Normal Form in presence of torsion
%%%%%%%%%%%%%%%%%%%%%%%%%%%%%%%%%%%%%%%%%%%%%%%%%%%%%%%%%%%%%%%%%%%%%%%%%%%%%%%%

%It is obvious that, if~$f$ has torsion~$\tau$, then we can always restrict ourselves to study the torsion-free~$\tau$-th iterate~$f^\tau$ of~$f$ {\sf non sono sicura che sia l'iterata $\tau$-esima quella da prendere perch\'e se esistono $m_1,\dots m_n$ interi non tutti nulli tali che
%$$
%m_1\phe_1+\cdots+m_n\phe_n = {1\over \tau}
%$$
%allora se moltiplico tutto per $\tau$ ottengo
%$$
%\tau m_1\phe_1+\cdots+\tau m_n\phe_n = {1},
%$$ 
%e allora $[\tau\phe]$ ha torsione $1$! L'iterata giusta dovrebbe essere quella che veniva fuori dalle $r$-uple ridotte, solo che c'\`e quell'esempio... bisogna capire che relazione c'\`e fra i vari numeri di queste $r$-uple particolari}, but in this case {\sf capire che diavolo significa da un punto di vista di forme normali... questo invece non dovrebbe essere infattibile, perch\'e se \`e giusto quello che ho scritto appena sopra, bisogna prendere l'iterata $m$-esima, e poi si ha che $[m\phe] =[\sum_{k=2}^r\beta_k m\eta^{(k)}]$; allora se per caso $f^m$ fosse normalizzabile, avremmo un'azione con pesi $\eta^{(2)},\dots, \eta^{(r)}$ e quindi, se $f^m$ commutava con quell'azione implica che deve commutarci anche $f$ vuol dire che avremo eliminato quasi tutti i pezzi che non ci interessavano (le risonanze di $f$ sono un sottoinsieme delle risonanze di $f^m$}

%However, we can say something also studying our germ~$f$ directly.

\sm Let us consider now~$[\phe]\in\C/\Z$, of toric degree~$1\le r\le n$ and let~$\theta^{(1)},\dots, \theta^{(r)}$ be a~$r$-tuple of toric vectors associated to~$[\phe]$ with toric coefficients~$\alpha_1,\dots,\alpha_r$ rationally dependent with $1$. We shall put
$$
\ca D(\alpha_1,\dots,\alpha_r)= \{M\in\Z^r\mid m_1\alpha_1+\cdots+m_r\alpha_r\in\Z\},
$$
and
$$
{\rm Adm}(\theta^{(1)}, \dots, \theta^{(r)})=\bigcup_{j=1}^n{\rm Adm}_j(\theta^{(1)}, \dots, \theta^{(r)}),
$$
where
$$
{\rm Adm}_j(\theta^{(1)}, \dots, \theta^{(r)})= \{M\in\Z^r\mid \exists Q\in\N^n, |Q|\ge 2~{\rm s.t.}~m_k = \langle Q-e_j,\theta^{(k)}\rangle\,\forall k=1, \dots, r\}\cup\{O\},
$$
for any~$j\in\{1,\dots, n\}$. 

\sm Even if, in this case, it is not always true that we can compute the resonances of $[\phe]$ as intersection of additive resonances, we can say many things on the resonant multi-indices using reduced $r$-tuples associated to $[\phe]$.

\sm\thm{Le4.3}{Lemma} {\sl Let~$[\phe]\in(\C/\Z)^n$ be of toric degree~$1\le r\le n$ and let~$\eta^{(1)},\dots, \eta^{(r)}$ be a reduced~$r$-tuple of toric vectors associated to~$[\phe]$ with toric coefficients~$1/m,\beta_2,\dots,\beta_r$. Then~{\parindent=30pt
\sm\item{(i)} $\ca D(1/m,\beta_2,\dots,\beta_r)= \{(hm,0, \dots, 0)\mid h\in\Z \}\subset\Z^r$;
\sm\item{(ii)} we have
$$
\ca D(1/m,\beta_2,\dots,\beta_r)\cap {\rm Adm}(\eta^{(1)}, \dots, \eta^{(r)})\ne \{O\}
$$
if and only there exist~$Q\in\N^n$, with~$|Q|\ge 2$ and $j\in\{1,\dots, n\}$ such that
$$
\la Q - e_j, \eta^{(1)}\ra \in m\Z\setminus\{0\} \quad{\rm and}\quad Q\in\bigcap_{k=2}^r \res^+_j(\eta^{(k)});
$$
\sm\item{(iii)} we have
$$
\res_j([\phe])=\{Q\in\N^n\mid |Q|\ge 2, \la Q - e_j, \eta^{(1)}\ra \in m\Z\}\cap \bigcap_{k=2}^r \res_j^+(\eta^{(k)}),
$$
for any $j\in\{1,\dots, n\}$.
In particular,
$$ 
\bigcap_{k=2}^r \res_j^+(\eta^{(k)}) \supseteq \res_j([\phe])\supseteq
\bigcap_{k=1}^r \res_j^+(\eta^{(k)}),\tag{trecasi}
$$
for any $j\in\{1,\dots, n\}$.
\sm\item{(iv)} $[\phe_j]=[\phe_h]$ implies that~$m$ divides $\eta^{(1)}_j-\eta^{(1)}_h$, and that~$\eta^{(k)}_j=\eta^{(k)}_h$ for any~$k=2, \dots, r$.\sm}}

\sm\proof (i) One inclusion is obvious. Conversely, let~$M\in \ca D(1/m,\beta_2,\dots,\beta_r)$; then
$$
m_1{1\over m}+ m_2\beta_2+\cdots+m_r\beta_r\in\Z.
$$
Since~$\beta_2, \dots, \beta_r$ are rationally independent with $1$, this implies~$m_2=\cdots= m_r=0$, thus we must have
%$$m_2\beta_2+\cdots+m_r\beta_r\ne \left[{1\over m}\right],$$
%for every~$(m_2, \dots, m_r)\in\Z^{r-1}$, so it must be
$m_1/m\in\Z$, and we are done.

\sm(ii) It is immediate from the definitions of~$\ca D(1/m,\beta_2,\dots,\beta_r)$ and~${\rm Adm}(\eta^{(1)}, \dots, \eta^{(r)})$ and from (i).

\sm(iii) It is immediate from (ii) and from 
$$
\left[\la Q,\phe \ra - \phe_j\right] = \left[{1\over m} \la Q-e_j,\eta^{(1)} \ra + \sum_{k=2}^r \beta_k \la Q-e_j,\eta^{(k)} \ra \right]. \tag{eq4.1}
$$ 

\sm(iv) If~$[\phe_j]=[\phe_h]$, then
$$
\left[{1\over m}\eta^{(1)}_j + \beta_2\eta^{(2)}_j + \cdots + \beta_r\eta^{(r)}_j \right]
=
\left[{1\over m}\eta^{(1)}_h + \beta_2\eta^{(2)}_h + \cdots + \beta_r\eta^{(r)}_h\right],
$$
hence
$$
{1\over m}\left(\eta^{(1)}_j-\eta^{(1)}_h \right) +\beta_2\left(\eta^{(2)}_j-\eta^{(2)}_h \right) \cdots + \beta_r\left(\eta^{(r)}_j-\eta^{(r)}_h\right)\in\Z,
$$
and, since~$\eta^{(k)}_j-\eta^{(k)}_h\in\Z$ for~$k=1,\dots, r$, the assertion follows as in (i).
\qed

\sm\thm{ReCompRidotta}{Remark} Note that, given~$[\phe]\in(\C/\Z)^n$ of toric degree~$1\le r\le n$, if~$\eta^{(1)},\dots, \eta^{(r)}$ is a reduced~$r$-tuple of toric vectors associated to~$[\phe]$ with toric coefficients~$1/m,\beta_2,\dots,\beta_r$, and such that~$[\phe_j]=[\phe_h]$ for some distinct coordinates $j$ and $h$, but $\eta^{(1)}_j\ne\eta^{(1)}_h$, then, since $m$ divides $\eta^{(1)}_j-\eta^{(1)}_h$, we have 
$$
{1\over m} \eta^{(1)}_j = {1\over m} \eta^{(1)}_h + {1\over m} \left(\eta^{(1)}_j-\eta^{(1)}_h\right);
$$
thus
$$
[\phe] = \left[{1\over m} \widetilde\eta^{(1)} + \sum_{k=2}^r\beta_k\eta^{(k)}\right]
$$
where, $\widetilde \eta^{(1)}_p= \eta^{(1)}_p$ for any $p\ne j,h$ and $\widetilde \eta^{(1)}_j=\widetilde \eta^{(1)}_h$, that is $\widetilde\eta^{(1)} = \eta^{(1)} - (\eta^{(1)}_j-\eta^{(1)}_h) e_j$, obtaining a compatible reduced $r$-tuple.

%{\bf Aggiungere Remark: Se ho una $r$-upla torica ridotta e ho $[\phe_j]=[\phe_h]$, ma $\eta^{(1)}_j\ne\eta^{(1)}_h$, allora, siccome~$m$ divide $\eta^{(1)}_j-\eta^{(1)}_h$, si ha 
%$$
%{1\over m} \eta^{(1)}_j = {1\over m} \eta^{(1)}_h + {1\over m} \left(\eta^{(1)}_j-\eta^{(1)}_h\right)
%$$
%e posso prendere un nuovo $\eta^{(1)}$ con $\eta^{(1)}_j= \eta^{(1)}_h$, ottenendo una $r$-upla compatibile.}

\sm Even in the torsion case, toric $r$-tuples associated to a same vector $[\phe]$ have to verify certain properties on the resonances, as next result shows.

\sm\thm{LeRapprRid}{Lemma} {\sl Let~$[\phe]\in(\C/\Z)^n$ be of toric degree~$1\le r\le n$ and in the torsion case. Let~$\eta^{(1)},\dots, \eta^{(r)}$ be a reduced~$r$-tuple of toric vectors associated to~$[\phe]$ with toric coefficients~$1/m,\beta_2,\dots,\beta_r$ and let~$\xi^{(1)},\dots, \xi^{(r)}$ be a reduced~$r$-tuple of toric vectors associated to~$[\phe]$ with toric coefficients~$1/\widetilde m,\gamma_2,\dots,\gamma_r$. Then we have
$$
\bigcap_{k=2}^r\res_j^+(\eta^{(k)}) = \bigcap_{k=2}^r\res_j^+(\xi^{(k)}),
$$
for any $j=1, \dots, n$.}

\sm\proof We have
$$
[\phe]
= \left[{1\over m}\eta^{(1)}+ \sum_{k=2}^r \beta_k\eta^{(k)}\right]
= \left[{1\over \widetilde m}\xi^{(1)}+ \sum_{k=2}^r \gamma_k\xi^{(k)}\right].
$$
Then
$$
[m\widetilde m\phe]
= \left[\sum_{k=2}^r m\widetilde m\beta_k\eta^{(k)}\right]
= \left[\sum_{k=2}^r m\widetilde m\gamma_k\xi^{(k)}\right],
$$
and, by Proposition \rf{PrCriterio}, $[m\widetilde m\phe]$ has toric degree $r-1$ and is torsion-free, because $\beta_2,\dots, \beta_r$ and $\gamma_2,\dots, \gamma_r$ are rationally independent with $1$. Therefore, by Lemma \rf{Le1.2}, we have
$$
\bigcap_{k=2}^r\res_j^+(\eta^{(k)}) = \res_j([m\widetilde m\phe]) = \bigcap_{k=2}^r\res_j^+(\xi^{(k)}),
$$
for any $j=1, \dots, n$, and we are done. \qed

\sm As Theorem \rf{PrNoImpureTorsionFree} shows, it is not possible that $\la P, \eta^{(1)}\ra \in m\Z$ for any $P\in\Z^n$ such that~$\la P,\eta^{(k)}\ra =0$ for $k=2, \dots, r$. However, it is possible that
%$$
%\ca D(1/m,\beta_2,\dots,\beta_r)\cap {\rm Adm}(\eta^{(1)}, \dots, \eta^{(r)})\ne \{O\},
%$$
%and
$$
\res_j([\phe]) = \bigcap_{k=2}^r\res_j^+(\eta^{(k)})
$$
for any $j\in\{1,\dots, n\}$, as next example shows.

\sm\thm{ExImpureTorsion2}{Example} Let us consider the vector
$$
[\phe]= \left[{1\over 3}\pmatrix{0 \cr
								 0 \cr
								1 \cr
						    1 } + {\sqrt{2}}\pmatrix{-12 \cr
						    			         0 \cr
									           0\cr
									           1} 
									+ {\sqrt{3}}\pmatrix{0 \cr
						    			        5 \cr
									          2 \cr
									          0}\right]\in(\C/\Z)^4,
$$
of toric degree $3$.
In this case $\ca D(1/3,\sqrt{2}, \sqrt{3}) = \{(3h, 0, 0)\mid h\in\Z\}$. 
We have
$$
\la P,\eta^{(2)}\ra = -12 p_1 + p_4= 0
$$
if and only if
$$
P\in \pmatrix{ 1 \cr
 			  0 \cr
			  0\cr
			  12}\Z\oplus e_2\Z\oplus e_3\Z,
$$
and
$$
\la P,\eta^{(3)}\ra = 5 p_2 + 2 p_3= 0
$$
if and only if
$$
P\in \pmatrix{ 0 \cr
 			  -2 \cr
			  5\cr
			  0}\Z\oplus e_1\Z\oplus e_4\Z.
$$
We have
$$
\eqalign{
&\res_1^+(\eta^{(2)}) =\left\{(q_1, q_2, q_3, 12 (q_1-1) )\mid q_1, q_2, q_3\in\N, 13q_1+q_2+q_3\ge 14\right\}\cr
&\res_2^+(\eta^{(2)}) =\left\{(q_1, q_2, q_3, 12 q_1 )\mid q_1, q_2, q_3\in\N, 13q_1+q_2+q_3\ge 2\right\}\cr
&\res_3^+(\eta^{(2)})= \res_2^+(\eta^{(2)})\cr
&\res_4^+(\eta^{(2)}) =\left\{(q_1, q_2, q_3, 12q_1 +1 )\mid q_1, q_2, q_3\in\N, 13q_1+q_2+q_3\ge 1\right\},}
$$
and
$$
\eqalign{
&\res_1^+(\eta^{(3)}) =\left\{(q_1, 0, 0, q_4 )\mid q_1, q_4\in\N, q_1+q_4\ge 2\right\}\cr
&\res_2^+(\eta^{(3)}) =\left\{(q_1, 1, 0, q_4 )\mid q_1, q_4\in\N, q_1+q_4\ge 1\right\}\cr
&\res_3^+(\eta^{(3)})=\left\{(q_1, 0, 1, q_4 )\mid q_1, q_4\in\N, q_1+q_4\ge 1\right\}\cr
&\res_4^+(\eta^{(3)}) =\res_1^+(\eta^{(3)}).}
$$
Moreover for each multi-index of the form $(p, 0, 0, 12p)$ with $p\ge 1$, we get
$$
\left\la \pmatrix{ p \cr
 			  0 \cr
			  0\cr
			  12p}, \eta^{(1)}\right\ra = 12p\in 3\,\Z.
$$
Then 
%$$
%\ca D (1/3, \sqrt{2},\sqrt{3}) \cap {\rm Adm}_j(\eta^{(1)}, \eta^{(2)},\eta^{(3)})\ne\{O\}
%$$
it is easy to verify that
$$
\res_j([\phe]) = \res_j^+(\eta^{(2)})\cap\res_j^+(\eta^{(3)}),
$$
for $j=1,\dots,4$.

%$$
%\eqalign{
%&\res_1([\phe])=\res_1^+(\eta^{(2)})\cap\res_1^+(\eta^{(3)}) =\left\{(q, 0,0, 12 (q-1) )\mid q\in\N, 13q\ge 1\right\}\cr
%&\res_2([\phe])=\res_2^+(\eta^{(2)})\cap\res_2^+(\eta^{(3)}) =\left\{(q, 1,0, 12 q )\mid q\in\N, 13q\ge 1\right\}\cr
%&\res_3([\phe])=\res_3^+(\eta^{(2)})\cap\res_3^+(\eta^{(3)}) =\left\{(q, 0,1, 12 q )\mid q\in\N, 13q\ge 1\right\}\cr
%&\res_4([\phe])=\res_4^+(\eta^{(2)})\cap\res_4^+(\eta^{(3)}) =\left\{(q, 0,0, 12 q + 1 )\mid q\in\N, 13q\ge 1\right\}.}
%$$

\sm\thm{ReEcallePerMarco}{Remark} Last example shows that, even in the torsion case, there are vectors $[\phe]\in(\C/\Z)^n$ such that, for any $j$, $\res_j([\phe])$ can be written as intersection of sets of additive resonances.

\sm We have then the following definition.

%\sm\defin{ImpureTorsion} Let~$[\phe]\in(\C/\Z)^n$ be of toric degree~$1\le r\le n$ in the torsion case. We say that $[\phe]$ is in the {\it impure torsion case} if, given~$\eta^{(1)},\dots, \eta^{(r)}$ a reduced~$r$-tuple of toric vectors associated to~$[\phe]$ with toric coefficients~$1/m,\beta_2,\dots,\beta_r$, we have
%$$
%\ca D(1/m,\beta_2,\dots,\beta_r)\cap {\rm Adm}_j(\eta^{(1)}, \dots, \eta^{(r)})= \{O\}.
%$$
%Otherwise we say that $[\phe]$ is in the {\it pure torsion case}.

%\sm\defin{ImpureTorsion} Let~$[\phe]\in(\C/\Z)^n$ be of toric degree~$1\le r\le n$ in the torsion case. Let~$\eta^{(1)},\dots, \eta^{(r)}$ be a reduced~$r$-tuple of toric vectors associated to~$[\phe]$ with toric coefficients~$1/m,\beta_2,\dots,\beta_r$.  We say that the toric $r$-tuple is {\parindent=30pt
%\sm\item{(i)} in the {\it impure torsion case of type $1$} if we have
%$$
%\ca D(1/m,\beta_2,\dots,\beta_r)\cap {\rm Adm}_j(\eta^{(1)}, \dots, \eta^{(r)})= \{O\};
%$$
%\sm\item{(ii)} in the {\it impure torsion case of type $2$} if we have
%$$
%\ca D(1/m,\beta_2,\dots,\beta_r)\cap {\rm Adm}_j(\eta^{(1)}, \dots, \eta^{(r)})={\rm Adm}_j(\eta^{(1)}, \dots, \eta^{(r)})\ne\{O\};
%$$
%\sm\item{(iii)} in the {\it pure torsion case} if we have
%$$
%\ca D(1/m,\beta_2,\dots,\beta_r)\cap {\rm Adm}_j(\eta^{(1)}, \dots, \eta^{(r)})\ne{\rm Adm}_j(\eta^{(1)}, \dots, \eta^{(r)})
%$$
%and
%$$
%\ca D(1/m,\beta_2,\dots,\beta_r)\cap {\rm Adm}_j(\eta^{(1)}, \dots, \eta^{(r)})\ne \{O\}.
%$$ 
%\sm}

%\sm Example \rf{ExImpureTorsion2} shows that the impure torsion case can occur.

\sm\defin{ImpureTorsion} Let~$[\phe]\in(\C/\Z)^n$ be of toric degree~$1\le r\le n$ and in the torsion case. We say that $[\phe]$ is in the {\it impure torsion case} if, given~$\eta^{(1)},\dots, \eta^{(r)}$ a reduced~$r$-tuple of toric vectors associated to~$[\phe]$ with toric coefficients~$1/m,\beta_2,\dots,\beta_r$, we have
$$
\res_j([\phe]) =\bigcap_{k=2}^r\res_j^+(\eta^{(k)}),\tag{torsioneimpura2}
$$
for any $j\in\{1,\dots, n\}$. Otherwise we say that $[\phe]$ is in the {\it pure torsion case}.

\sm The next result shows that the impure torsion case is well-defined, i.e., it does not depend on the chosen toric $r$-tuple.

\sm\thm{BuonaDefTorsioneType2}{Lemma} {\sl Let~$[\phe]\in(\C/\Z)^n$ be of toric degree~$1\le r\le n$ and in the torsion case. Let~$\eta^{(1)},\dots, \eta^{(r)}$ be a reduced~$r$-tuple of toric vectors associated to~$[\phe]$ with toric coefficients~$1/m,\beta_2,\dots,\beta_r$. If  
$$
\res_j([\phe]) =\bigcap_{k=2}^r\res_j^+(\eta^{(k)}),\tag{torsioneimpura2}
$$
for any $j\in\{1,\dots, n\}$, then \rf{torsioneimpura2} holds for any other reduced toric $r$-tuple associated to $[\phe]$.}

\sm\proof Let~$\xi^{(1)},\dots, \xi^{(r)}$ be another reduced~$r$-tuple of toric vectors associated to~$[\phe]$ with toric coefficients~$1/\widetilde m,\gamma_2,\dots,\gamma_r$. 
Since $\eta^{(1)},\dots, \eta^{(r)}$ is in the impure torsion case, we have
$$
\res_j([\phe]) = \bigcap_{k=2}^r \res_j^+(\eta^{(k)}),
$$
but, thanks to Lemma \rf{LeRapprRid}, we have
$$
\bigcap_{k=2}^r\res_j^+(\eta^{(k)}) = \bigcap_{k=2}^r\res_j^+(\xi^{(k)}),
$$
for any $j=1, \dots, n$, that is $\xi^{(1)},\dots, \xi^{(r)}$ satisfy \rf{torsioneimpura2}. \qed

\sm\defin{ImpureGerm} Let~$f$ be a germ of biholomorphism of~$\C^n$ fixing the origin. We say that $f$ is in the {\it impure torsion case} [resp., {\it in the pure torsion case}] if, denoting with $\lambda=\{\lambda_1, \dots, \lambda_n\}$ the spectrum of~$\d f_O$, the unique $[\phe]\in(\C/\Z)^n$ such that~$\lambda=e^{2\pi i [\phe]}$ is in the impure torsion case [resp., in the pure torsion case].

\sm\thm{Pr4.1}{Theorem} {\sl Let $f$ be a germ of biholomorphism of~$\C^n$ fixing the origin~$O$ of toric degree~$1\le r\le n$ and in the impure torsion case. Then it admits a holomorphic Poincar\'e-Dulac normalization if and only if there exists a holomorphic effective action on~$(\C^n,O)$ of a torus of dimension~$r-1$ commuting with~$f$, and such that the columns of the weight matrix of the action are reduced torsion-free toric vectors associated to~$f$.}

\sm\proof It follows from Theorem \rf{Te2.1}, Lemma \rf{Le4.3} and Lemma \rf{BuonaDefTorsioneType2}. \qed 

% Nota Bene: La compatibilita' e' garantita dal fatto che usiamo solo i vettori che corrispondono ai beta_k e allora posso usare il lemma \rf{Le4.3}

The next examples show that, in case of pure torsion there are more possible cases.
 
\sm\thm{ExImpureTorsionFree1}{Example} Let us consider the vector
$$
[\phe]= \left[{1\over 6}\pmatrix{1 \cr
		   3 } + \sqrt{2}\pmatrix{1 \cr
		   6}
		   \right]\in(\C/\Z)^2,
$$
of toric degree $2$.
In this case $\ca D(1/6,\sqrt{2})= \{(6h, 0)\mid h\in\Z\}$. 
We have
$$
\la P,\eta^{(2)}\ra = p_1 + 6 p_2 = 0
$$
if and only if
$$
P\in \pmatrix{-6\cr
				1}\Z,
$$
hence
$$
\res_1^+(\eta^{(2)}) = \void \quad\hbox{and}\quad \res_2^+(\eta^{(2)}) = \{(6, 0)\}.
$$
Since
$$
\la (6,-1), \eta^{(1)}\ra = 3 \not\in 6\Z,
$$
we have
$$
\ca D (1/6, \sqrt{2}) \cap {\rm Adm}_j(\eta^{(1)}, \eta^{(2)})= \{O\}
$$
for $j=1,2$, so we have
$$
\res_j([\phe]) =\bigcap_{k=1}^r\res_j^+(\eta^{(k)})=\void,
$$
for $j=1,2$. Moreover, it is evident that the torsion is $2$. %{\sf e la torsione e' 2}

\sm\thm{ExPurePresenceTorsion2}{Example} Let us consider the vector
$$
[\phe]= \left[{1\over 7}\pmatrix{1 \cr
		   3 } + \sqrt{2}\pmatrix{1 \cr
		   -6}
		   \right]\in(\C/\Z)^2,
$$
of toric degree $2$.
In this case $\ca D(1/7,\sqrt{2}) = \{(7h, 0)\mid h\in\Z\}$. 
We have
$$
\res_1^+(\eta^{(1)}) = \void \quad\hbox{and}\quad \res_2^+(\eta^{(1)}) = \{(3,0)\},
$$
and
$$
\res_1^+(\eta^{(2)}) = \{(6h+1,h)\mid h\in\N\setminus\{0\}\} \quad\hbox{and}\quad \res_2^+(\eta^{(2)}) = \{(6h, h+1)\mid h\in\N\setminus\{0\}\};
$$
then
$$
\res_1^+(\eta^{(1)})\cap \res_1^+(\eta^{(2)})=\void \quad\hbox{and}\quad \res_2^+(\eta^{(1)})\cap \res_2^+(\eta^{(2)})=\void.
$$
However, we have
$$
\la (6h,h), \eta^{(1)}\ra \in 9\Z;
$$
hence
%$$
%\ca D (1/7, \sqrt{2}) \cap {\rm Adm}_j(\eta^{(1)}, \eta^{(2)})\ne \{O\}
%$$
%for $j=1,2$, we are in the pure torsion case, and 
we have
$$
\eqalign{
&\res_1^+(\eta^{(2)})\supset\res_1([\phe]) =\left\{(42 h +1, 7 h)\mid h\in\N\setminus\{0\}\right\}\supset \res_1^+(\eta^{(1)})\cap \res_1^+(\eta^{(2)})\cr
&\res_2^+(\eta^{(2)})\supset\res_2([\phe]) = \left\{(42 h, 7 h +1)\mid h\in\N\setminus\{0\}\right\}\supset \res_2^+(\eta^{(1)})\cap \res_2^+(\eta^{(2)}).}
$$
Moreover, it is not difficult to verify that the torsion is $7$.
%{\sf e la torsione e' $7$.}

\sm In the pure torsion case, one could ask whether, given a toric $r$-tuple $\eta^{(1)},\dots, \eta^{(r)}$ associated to $[\phe]$ such that
$$
\bigcap_{k=2}^r\res_j^+(\eta^{(k)})\supset\res_j([\phe]) \supset \bigcap_{k=1}^r\res_j^+(\eta^{(k)}),\tag{torsionepura}
$$
for some $j\in\{1,\dots, n\}$, then this is true for any other toric $r$-tuple associated to $[\phe]$. This is not always true, as next example shows.

\sm\thm{ExPureTorsion3}{Example} Let us consider the vector
$$
[\phe]= \left[{1\over 3}\pmatrix{1 \cr
								 1 \cr
								 1 \cr
								 1}
		  + \sqrt{2}\pmatrix{1 \cr
		                     6 \cr
		                     0 \cr
		                     0}
		   + \sqrt{3}\pmatrix{0 \cr
		                     0 \cr
		                     -1 \cr
		                     5}
		   \right]\in(\C/\Z)^4,
$$
of toric degree $3$.
In this case $\ca D(1/3,\sqrt{2},\sqrt{3}) = \{(3h, 0,0)\mid h\in\Z\}$. 
We have
$$
\res_j^+(\eta^{(1)}) =\void,
$$
for $j=1,\dots, 4$,
$$
\eqalign{
&\res_1^+(\eta^{(2)}) =\left\{(1,0,p,q)\mid p,q\in\N, p + q\ge 1\right\}\cr
&\res_2^+(\eta^{(2)}) =\left\{(6,0,p, q)\mid p,q\in\N\right\}\cup\left\{(0,1,p, q)\mid p,q\in\N, p+q\ge 1\right\} \cr
&\res_3^+(\eta^{(2)})= \left\{(0,0,p, q)\mid p,q\in\N, p + q \ge 2\right\}\cr
&\res_4^+(\eta^{(2)}) =\res_3^+(\eta^{(2)}),}
$$
and
$$
\eqalign{
&\res_1^+(\eta^{(3)}) =\left\{(h, k, 5q,q)\mid h,k, q\in\N, h+k + 6q\ge 2\right\}\cr
&\res_2^+(\eta^{(3)}) =\res_1^+(\eta^{(3)})\cr
&\res_3^+(\eta^{(3)})= \left\{(h,k,5q + 1, q)\mid h,k,q\in\N, h+k+ 6q \ge 1\right\}\cr
&\res_4^+(\eta^{(3)})= \left\{(h,k,5(q -1), q)\mid h,k,q\in\N, h+k+ 6q \ge 7\right\}.}
$$
Then we have
$$
\bigcap_{k=1}^3\res_j^+(\eta^{(k)}) = \void,
$$
for $j=1,\dots, 4$, but it is not difficult to verify that
$$
\res_2([\phe])= \left\{(0,1,5q, q)\mid q\in\N^*\right\}\ne\void\quad\hbox{and}\quad\res_j([\phe])= \res_j^+(\eta^{(2)})\cap\res_j^+(\eta^{(3)})~~j=1,3,4.
$$
Then, since
$$
\res_2^+(\eta^{(2)})\cap \res_2^+(\eta^{(3)}) = \left\{(6, 0, 5q,q)\mid q\in\N\right\} \cup \left\{(0,1,5q, q)\mid q\in\N^*\right\} \ne\res_2([\phe]),
$$
we are in the pure torsion case, but we cannot write all the resonances of $[\phe]$ as intersection of the additive resonances of $\eta^{(1)}$, $\eta^{(2)}$ and $\eta^{(3)}$.
However, we can write
$$
[\phe]= \left[{1\over 3}\pmatrix{1 \cr
								 -2\cr
								 1 \cr
								 -5}
		  + \sqrt{2}\pmatrix{1 \cr
		                     6 \cr
		                     0 \cr
		                     0}
		   + \sqrt{3}\pmatrix{0 \cr
		                     0 \cr
		                     -1 \cr
		                     5}
		   \right],
$$
and it is not difficult to verify that, in this representation, we have
$$
\res_j([\phe])= \bigcap_{k=1}^3\res_j^+(\xi^{(k)}),
$$
for $j=1,\dots, 4$.

\sm\thm{ReTorsionePuraCasino}{Example} If $[\phe]\in(\C/\Z)^2$ is given by Example \rf{ExPurePresenceTorsion2}, we saw that we can write it in the form
$$
[\phe] = \left[{1\over \tau}\eta^{(1)} +\beta\eta^{(2)}\right]
$$
so that
$$
\res_j^+(\eta^{(2)})\supset\res_j([\phe]) \supset \res_j^+(\eta^{(1)})\cap \res_j^+(\eta^{(2)}),\tag{torsionepura3}
$$
for all $j$.  Furthermore, it is easy to check that $[\phe]$ does not admit any
reduced representation
$$
[\phe] = \left[{1\over \tau q}\xi^{(1)} +\gamma\xi^{(2)}\right]
$$
such that for any $j$ we have
$$
\res_j([\phe])= \res_j^+(\xi^{(1)})\cap  \res_j^+(\xi^{(2)}).\tag{torsionepura4}
$$

\sm We are then led to the following

\sm\defin{ReduciblePureTorsion} Let~$[\phe]\in(\C/\Z)^n$ be of toric degree~$1\le r\le n$ and in the pure torsion case. We say that $[\phe]$ {\it can be simplified} if it admits a reduced~$r$-tuple of toric vectors~$\eta^{(1)},\dots, \eta^{(r)}$ such that
$$
\res_j([\phe]) =\bigcap_{k=1}^r\res_j^+(\eta^{(k)}),\tag{torsionepura3}
$$
for all $j=1,\dots, n$. The $r$-tuple $\eta^{(1)},\dots, \eta^{(r)}$ is said a {\it simple reduced} $r$-tuple associated to $[\phe]$.

\sm\defin{ReduciblePureGerm} Let~$f$ be a germ of biholomorphism of~$\C^n$ fixing the origin in the pure torsion case. We say that $f$ {\it can be simplified}   if, denoting with $\lambda=\{\lambda_1, \dots, \lambda_n\}$ the spectrum of~$\d f_O$, the unique $[\phe]\in(\C/\Z)^n$ such that~$\lambda=e^{2\pi i [\phe]}$ can be simplified.

%\sm\thm{Pr4.1}{Theorem} {\sl Let $f$ be a germ of biholomorphism of~$\C^n$ fixing the origin~$O$ of toric degree $1\le r\le n$ in the impure torsion case. Then it admits a holomorphic Poincar\'e-Dulac normalization if and only if there exists a holomorphic effective action on~$(\C^n,O)$ of a torus of dimension~$r$ commuting with~$f$ and such that the columns of the weight matrix of the action are a reduced~$r$-tuple of toric vectors associated to~$f$.}

%\sm\proof It follows from Theorem \rf{Te2.1}, Lemma \rf{Le4.3} and Remark \rf{Re7.1}. \qed 

\sm\thm{Pr4.2}{Theorem} {\sl Let $f$ be a germ of biholomorphism of~$\C^n$ fixing the origin~$O$ of toric degree~$1\le r\le n$ and in the pure torsion case and it can be simplified. Then:{\parindent=30pt
\sm\item{(i)} if~$\d f_O$ is diagonalizable, $f$ admits a holomorphic Poincar\'e-Dulac normalization if and only if there exists a holomorphic effective action on~$(\C^n,O)$ of a torus of dimension~$r$ commuting with~$f$ and such that the columns of the weight matrix $\Theta$ of the action are a simple reduced~$r$-tuple of toric vectors associated to~$f$;
\sm\item{(ii)} if~$\d f_O$ is not diagonalizable and there exists a simple reduced~$r$-tuple of toric vectors associated to~$f$ such that its vectors are the columns of a matrix $\Theta$ compatible with~$\d f_O$,~$f$ admits a holomorphic Poincar\'e-Dulac normalization if and only if there exists a holomorphic effective action on~$(\C^n,O)$ of a torus of dimension~$r$ commuting with~$f$ and with weight matrix $\Theta$.
\sm}}

\sm\proof It follows from Theorem \rf{Te2.1}. \qed 

\sm\thm{ReComp1}{Remark} Note that we cannot get rid of the compatibility hypothesis in the case of $\d f_O$ non diagonalizable, because if we change a simple reduced toric $r$-tuple as in Remark \rf{ReCompRidotta}, it is not true that we obtain another simple reduced $r$-tuple. In fact, if~$[\phe]\in(\C/\Z)^n$ has toric degree~$1\le r\le n$, and~$\eta^{(1)},\dots, \eta^{(r)}$ is a simple reduced~$r$-tuple of toric vectors associated to~$[\phe]$ with toric coefficients~$1/m,\beta_2,\dots,\beta_r$, but we have~$[\phe_j]=[\phe_h]$ for some distinct coordinates $j$ and $h$, and $\eta^{(1)}_j\ne\eta^{(1)}_h$, then for every $P\in\res_l([\phe])$, the equality
$$
[\phe] = \left[{1\over m} \widetilde\eta^{(1)} + \sum_{k=2}^r\beta_k\eta^{(k)}\right]
$$
with $\widetilde\eta^{(1)} = \eta^{(1)} - (\eta^{(1)}_j-\eta^{(1)}_h) e_j$, only implies
$$
{\eta^{(1)}_j-\eta^{(1)}_h\over m} (\delta_{lh}-p_h)\in \Z
$$
and there well can be resonant multi-indices with $p_h\ne 1$. 

%\me Let us consider now the case of~$[\phe]\in(\C/\Z)^n$ be of toric degree~$1\le r\le n$ having a~$r$-tuple of toric vectors~$\theta^{(1)},\dots, \theta^{(r)}$ associated to~$[\phe]$ with toric coefficients~$\alpha_1,\dots,\alpha_r$ rationally dependent with $1$ and~$\res_j([\phe]) \supset \bigcap_{k=1}^r \res_j^+(\theta^{(k)})$ for some~$j\in\{1,\dots, n\}$, i.e., thanks to Lemma \rf{Le4.2}, 
%$$\ca D(\alpha_1,\dots,\alpha_r)\cap {\rm Adm}_j(\theta^{(1)}, \dots, \theta^{(r)})\ne \{O\}.$$
%This case is the worst one, but using Lemma \rf{Pr1.0}, we can say something on the resonant multi-indices.

%However, as we already observed in Remark \rf{Re7.1}, in the impure torsion case, even if we have \'Ecalle's torsion we can compute the resonances of $[\phe]$ which are multiplicative, as intersection of additive resonances. {\sf inserire gli esempi dei casi del paragrafo 2}

%\sm\thm{ExTorsion}{Example} The vector in Example \rf{ExImpureTorsionFree1} has torsion $2$, but, since it is in the impure torsion case, we can compute its resonances as intersection of additive ones. 

%\sm\thm{ExTorsion2}{Example} The vector in Example \rf{ExPurePresenceTorsion2} has torsion $7$, and we are in the torsion case. 

\sm In case of pure torsion that cannot be simplified, we have the following result.

\sm\thm{PrPureTorsion}{Proposition} {\sl Let $f$ be a germ of biholomorphism of~$\C^n$ fixing the origin~$O$ of toric degree $1\le r\le n$ and in the pure torsion case and it cannot be simplified. If there exists a holomorphic effective action on~$(\C^n,O)$ of a torus of dimension~$r$ commuting with~$f$ and such that the columns of the weight matrix of the action are a reduced~$r$-tuple of toric vectors associated to~$f$, then $f$ admits a holomorphic Poincar\'e-Dulac normalization.}

\sm\proof It follows from Theorem \rf{Te2.1} and Lemma \rf{Le4.3}. \qed 

%\sm We only have a sufficient condition for holomorphic normalization, which is not necessary as next Example shows. 

%\sm\thm{EsNonSimpl}{Example} Let us consider $f$ the germ of biholomorphism of $(\C^2,O)$ given by
%$$
%\eqalign{&f_1(z,w)=\lambda z + z^{43}w^7 \cr
%  		 &f_2(z,w)=\mu w}
%$$
%where 
%$$
%\lambda = e^{2\pi i [{1\over 7} + \sqrt{2}]} \quad\hbox{and}\quad \mu = e^{2\pi i [{3\over 7} -6 \sqrt{2}]}.
%$$
%It is not difficult to verify that $f$ is in Poincar\'e-Dulac normal form, but since
%$$
%[\phe]= \left[{1\over 7}\pmatrix{1 \cr
%		   3 } + \sqrt{2}\pmatrix{1 \cr
%		   -6}
%		   \right]\in(\C/\Z)^2,
%$$
%is not simplifible (see Example \rf{ReTorsionePuraCasino}), $f$ cannot commute with any linear action of torus of dimension $2$ such that the columns of the weight matrix are a reduced couple of toric vectors associated to~$f$.

%Rimane la rottura del lineare, cioe' potrebbe commutare con una che va bene (e allora certo la forma normale di poincar\'e-dulac non e' piu' quella, anche se $f$ era gia' in forma di Poinca\'e-Dulac, hai capito?)

\me We end this section describing an algorithm to decide when a vector $[\phe]$ can be simplified.

\sm %Now, if we have torsion, we shall show in the next section that there exists a minimal $\widetilde m$ such that the $\widetilde m$-th iterate $f^{\widetilde m}$ of $f$ is torsion-free. 

We want to know when, given $[\phe]$ in the torsion case,
$$
[\phe] = \left[{1\over \tau p}\eta^{(1)} +\sum_{k=2}^r\beta_k\eta^{(k)}\right]
$$
of toric degree $r$, torsion $\tau\ge 2$, and such that there is $j\in\{1,\dots,n\}$ so that
$$
\bigcap_{k=1}^r\res_j^+(\eta^{(k)})\subset\res_j([\phe]) \subset \bigcap_{k=2}^r\res_j^+(\eta^{(k)}),\tag{torsionepura3}
$$
there is another reduced representation 
$$
[\phe] = \left[{1\over \tau q}\xi^{(1)} +\sum_{k=2}^r\gamma_k\xi^{(k)}\right]
$$
such that for any $j=1,\dots, n$ we have
$$
\res_j([\phe])= \bigcap_{k=1}^r\res_j^+(\xi^{(k)}).\tag{torsionepura4}
$$
We know that there must be $H\in\Z^n\setminus\{O\}$ such that
$$
{1\over \tau p}\eta^{(1)} +\sum_{k=2}^r\beta_k\eta^{(k)}= {1\over \tau q}\xi^{(1)} +\sum_{k=2}^r\gamma_k\xi^{(k)} + H.
$$
Since
$$
\bigcap_{k=2}^r\res_j^+(\eta^{(k)})= \bigcap_{k=2}^r\res_j^+(\xi^{(k)}),
$$
for any $j=1,\dots, n$, we have that 
$$
{1\over\tau p}\la\eta^{(1)},P-e_j\ra = {1\over\tau q}\la\xi^{(1)},P-e_j\ra + \la H,P-e_j\ra
$$
for any $P\in\bigcap_{k=2}^r\res_j^+(\eta^{(k)})$. Now, if $\la\xi^{(1)},P-e_j\ra =0$ it must be $\la\eta^{(1)},P-e_j\ra\in\tau p\Z$. On the contrary, if $\la\eta^{(1)},P-e_j\ra\in\tau p\Z$, then we would like to find $H$ such that $\la\xi^{(1)},P-e_j\ra =0$ that is, for any $j=1,\dots, n$,
$$
{1\over\tau p}\la\eta^{(1)},P-e_j\ra = \la H,P-e_j\ra
$$
for any $P\in\bigcap_{k=2}^r\res_j^+(\eta^{(k)})$ with $\la\eta^{(1)},P-e_j\ra\in\tau p\Z$. In fact, if such a vector exists, then, setting $q=p$, $\xi^{(1)} = \eta^{(1)}-\tau p H$, $\gamma_k=\beta_k$ and $\eta^{(k)} = \xi^{(k)}$ for $k=2,\dots, r$, we get
$$
[\phe]= \left[{1\over \tau p} \xi^{(1)} + \sum_{k=2}^r \gamma_k\xi^{(k)}\right],
$$
and for any $P\in\res_j([\phe])$ we have $P\in\bigcap_{k=2}^r\res_j^+(\xi^{(k)})$, and
$$
\la\xi^{(1)},P-e_j\ra = \la\eta^{(1)},P-e_j\ra - \la H,P-e_j\ra = 0,
$$
that is \rf{torsionepura4}.

\me We then have to study the structure of the intersection of a submodule of $\Z^n$ with $\N^n$. It turns out that such a structure is the following. We thank Jean \'Ecalle for suggesting the gist of the following argument.

%{\bf A occhio per capire se un vettore e' ``semplificabile'' o no si fa cosi'.
%Considero, per ogni $j=1, \dots, n$, l'insieme 
%$$
%\{P\in\Z^n\mid \la P-e_j,\eta^{(k)}\ra = 0~~\hbox{per}~~k=2,\dots, r\}= {\rm Span}_\Z\{P_{1,j},\dots, P_{n-r+1,j}\}
%$$
%e lo interseco con $\{Q\in\N^n\mid |Q|\ge 2\}$. Ho fatto un po' di prove e, anche se non ho una dimostrazione, dovrebbe essere vero che 
%$$
%\eqalign{
%\bigcap_{k=2}^r\res_j^+(\eta^{(k)})
%&={\rm Span}_\Z\{P_{1,j},\dots, P_{n-r+1,j}\}\cap \{Q\in\N^n\mid |Q|\ge 2\} \cr
%&= {\rm Span}_\N\{\widetilde P_{1,j},\dots, \widetilde P_{h(j),j}\}\cup\{Q_{1,j}, \dots, Q_{s(j),j}\}.}
%$$
%A questo punto considero l'insieme
%$$
%T_j=\{Q\in\N^n\mid \la Q-e_j, \eta^{(1)}\ra\in \tau\Z, |Q|\ge 2\}
%$$
%e ho
%$$
%\eqalign{
%T_j\cap\bigcap_{k=2}^r\res_j^+(\eta^{(k)})
%&=T_j\cap\left({\rm Span}_\N\{\widetilde P_{1,j},\dots, \widetilde P_{h(j),j}\}\cup\{Q_{1,j}, \dots, Q_{s(j),j}\}\right)\cr
%&={\rm Span}_\N\{P'_{1,j},\dots, P'_{t(j),j}\}\cup\{Q'_{1,j}, \dots, Q'_{u(j),j}\}.}
%$$
%e quindi $H$ dovra' verificare
%$$
%\la P'_{i,j}-e_j, H\ra = {1\over \tau} \la P'_{i,j}-e_j, \eta^{(1)}\ra
%$$ 
%per $i=1,\dots, t(j)$ e $j=1, \dots, n$ e 
%$$
%\la Q'_{i,j}-e_j, H\ra = {1\over \tau} \la Q'_{i,j}-e_j, \eta^{(1)}\ra
%$$
%per $i=1,\dots, u(j)$ e $j=1, \dots, n$.}

\sm Let ${\cal A}\subset\Z^n$ be a sub-module of $\Z^n$ where $n\in\N^*$, and let us denote by~$\ca A^+$ the set~$\ca A\cap \N^n$. 
For any vector $A=(a_1,\dots, a_n)\in\ca A$,  we denote by
$$
{\rm red}(A)= {1\over \alpha} A = \left({a_1\over \alpha},\dots, {a_n\over \alpha}\right)\tag{ecalle1}
$$
where $\alpha$ is the greatest common divisor of $a_1,\dots, a_n$. The {\it support} of a vector $A\in\Z^n$ is the set
$$
{\rm supp}(A)= \{j\in \{1,\dots, n\}\mid a_j\ne 0\} \subseteq \{1,\dots, n\}.
$$
Using the support we can then define a partial order on~$\ca A^+$ as follows: we say that $A\subseteq B$ if~${\rm supp}(A)\subset {\rm supp}(B)$, or the supports are equal and $A\le B$ in the usual lexicographic order.

\sm\defin{DeEcalle0} Let~$\ca A\subset\Z^n$ be any sub-module of $\Z^n$, where $n\in\N^*$, and let~$\ca A^+$ be the set~$\ca A\cap \N^n$. For any $A,B\in \ca A^+$ we define
$$
A/B = {\rm red}(q A-p B)\tag{ecalle2}
$$
where
$$
{p\over q} = \min_{j\in{\rm supp}(B)} \left({a_j\over b_j}\right).
$$
Obviously, if ${\rm supp}(B)\subseteq {\rm supp}(A)$, then $A/B \in\ca A^+$ and $A/B\subseteq A$.

\sm\defin{De1Ecalle} Let~$\ca A\subset\Z^n$ be any sub-module of $\Z^n$, where $n\in\N^*$, and let~$\ca A^+$ be the set~$\ca A\cap \N^n$. An element $M$ of~$\ca A^+$ is said {\it minimal} if it is minimal with respect to the partial order $\subseteq$.
An element $C$ of~$\ca A^+$ is said {\it cominimal} if for any minimal element $M$ of~$\ca A^+$ we have~$C-M\not \in\ca A^+$.

\sm Minimal elements have to satisfy certain properties.

\sm\thm{LeEcalle1}{Lemma} {\sl Let~$\ca A\subset\Z^n$ be any sub-module of $\Z^n$, where $n\in\N^*$, and let~$\ca A^+$ be the set~$\ca A\cap \N^n$. Two minimal elements of $\ca A^+$ have distinct supports.}

\sm\proof Let $M$ and $P$ be two distinct minimal elements of~$\ca A^+$ and suppose by contradiction that~${\rm supp}(M)= {\rm supp}(P)$. Then $A = M/P$ and $B=P/M$ both have supports strictly contained in the ones of $M$ and $P$ contradicting their minimality with respect to $\subseteq$.
\qed

\sm\thm{CoEcalle1}{Corollary} {\sl Let~$\ca A\subset\Z^n$ be any sub-module of $\Z^n$, where $n\in\N^*$, and let~$\ca A^+$ be the set~$\ca A\cap \N^n$. Then $\ca A^+$ contains only a finite number of minimal elements.}

\sm Minimal elements are a sort of generators of~$\ca A^+$ in a sense that next result clarifies.

\sm\thm{LeEcalle2}{Lemma} {\sl Let~$\ca A\subset\Z^n$ be any sub-module of $\Z^n$, where $n\in\N^*$, and let~$\ca A^+$ be the set~$\ca A\cap \N^n$. Then every element $A$ of~$\ca A^+$ can be written in the form
$$
A= {1\over \delta} (\alpha_1 M_1+\cdots+ \alpha_d M_d)\tag{ecalle5}
$$
where $\alpha_1, \dots, \alpha_d\in\N$, $M_1,\dots, M_d$ are the minimal elements, and $\delta=\delta(\ca A^+)\in\N^*$ depends only on~$\ca A^+$.}

\sm\proof If $A$ is non minimal, then there exists a minimal element $M_{j_1}\subseteq A$, and there exist~$\gamma_1,\delta_1\in\Q^+$ such that
$$
A = \gamma_1 M_{j_1} + \delta_1 A_1,
$$
where
$$
A_1 = A/M_{j_1},
$$
and ${\rm supp}(A_1)\subset {\rm supp}(A)$. If $A_1$ is not minimal, we can iterate this procedure getting
$$
A_1 = \gamma_2 M_{j_2} + \delta_2 A_2,
$$
with ${\rm supp}(A_2)\subset {\rm supp}(A_1)\subset {\rm supp}(A)$. The chain ${\rm supp}(A)\supset{\rm supp}(A_1)\supset {\rm supp}(A_3)\supset \cdots$ has to end because~$\ca A^+\subset \N^n$, then we eventually arrive to a decomposition of the form \rf{ecalle5}. Now~$\delta=\delta(\ca A^+)$ cannot be greater than the least common multiple of all $|\det({\rm M}^*)|$ where~${\rm M}^*$ varies in the square submatrices of order equal to the rank of the matrix having as columns all the minimal elements $M_1, \dots, M_d$ of $\ca A^+$.\qed

\sm The cominimal elements are finite too.

\sm\thm{LeEcalle3}{Lemma}  {\sl Let~$\ca A\subset\Z^n$ be any sub-module of $\Z^n$, where $n\in\N^*$, and let~$\ca A^+$ be the set~$\ca A\cap \N^n$. Then $\ca A^+$ contains only a finite number of cominimal elements.}

\sm\proof Let us assume by contradiction that there is an infinite sequence of distinct cominimal elements $\{C_j\}$. Thanks to Lemma \rf{LeEcalle2}, for each $j\ge 1$, we have
$$
C_j= {1\over \delta}\sum_{k=1}^d \gamma_{jk} M_k
$$
where $\gamma_{jk}\in\N$. Then there is an infinite subsequence $\{C_{j'}\}$ such that all the corresponding~$(\gamma_{j', 1},\dots, \gamma_{j',d})$ belong the a same class $(\gamma_1^*,\dots, \gamma_d^*)$ modulo $\delta\Z^d$.
Hence there is an infinite subsequence $\{C_{j''}\}$ such that at least one component $\gamma_{j'', k_0}$ diverges as  $j''$ tends to infinity, and such that the other components  $\gamma_{j'', k}$ with $k\ne k_0$ do not decrease. Then there exist at least two cominimal elements $C_{j_1}\le C_{j_2}$  such that
$$
C_{j_2} - C_{j_1} = \sum_{k=1}^d \tilde\gamma_k M_k
$$
with
$$
\tilde\gamma_k = {1\over \delta} \left(\gamma_{j_2,k} - \gamma_{j_1,k}\right)\in\N
$$
meaning that $C_{j_2}$ is not cominimal against the assumption. \qed

\sm For each element of $\ca A^+$, we want to find a decomposition with natural coefficients into linear combination of a finite number of elements of $\ca A^+$. This is possible using minimal and cominimal elements, as shown in next result.

\sm\thm{PrMinComin}{Proposition} {\sl Let~$\ca A\subset\Z^n$ be any sub-module of $\Z^n$, where $n\in\N^*$, and let~$\ca A^+$ be the set~$\ca A\cap \N^n$. Then for any $A\in\ca A^+$ there exist $l_1,\dots,l_d\in\N$ such that 
$$
A=\sum_{j=1}^d l_j M_j\tag{ecalle3}
$$
or 
$$
A=C_h+\sum_{j=1}^d l_j M_j\tag{ecalle4}
$$
for some $h\in\{1,\dots, e\}$, where $M_1,\dots, M_d$ are the minimal elements of~$\ca A^+$, and~$\{C_1,\dots, C_e\}$ are the cominimal elements of $\ca A^+$.}

\sm\proof If $A$ is non cominimal, there exists a minimal element $M_{j_1}\le A$; thus if~$A_1= A- M_{j_1}$ is not cominimal, we iterate the procedure. The chain $A\ge A_1\ge A_3\ge \cdots$ has to end with a zero, i.e., we get a decomposition of the form \rf{ecalle3}, or with a cominimal element, i.e., we get a decomposition of the form \rf{ecalle4}.\qed

\sm\thm{ReEcalle}{Remark} Note that it can happen that the number of minimal elements of $\ca A^+$ is not equal to the maximum number of $\Q$-linearly independent elements of $\ca A^+$. In fact, if we consider the submodule $\ca A$ of $\Z^4$ orthogonal to $(1,-1,-1,1)^T$, and $\ca A^+$, such a maximum is clearly $3$, but we have four minimal elements 
$$
\pmatrix{1\cr
		 1\cr
		 0\cr
		 0},
\pmatrix{0\cr
		 1\cr
		 0\cr
		 1},
\pmatrix{1\cr
		 0\cr
		 1\cr
		 0},
\pmatrix{0\cr
		 0\cr
		 1\cr
		 1},
$$
and we need all of them (and no cominimal) to ensure \rf{ecalle3} and \rf{ecalle4}.

\me Returning to our problem, if now we consider
$$
\ca A= \{Q\in\ \Z^n\mid \la Q, \eta^{(k)}\ra = 0,~\hbox{for}~k=2,\dots, r\},
$$
it is easy to verify that 
$$
\bigcap_{k=2}^r \res_j^+(\eta^{(k)})= \ca B^+_0 \cup \ca B^+_j
$$
where
$$
\ca B^+_0 = \{P\in\N^n\mid P=Q+e_j, Q\in\ca A^+, |Q|\ge 1\}
$$
and
$$
\ca B^+_j= \{P\in\N^n\mid P=Q+e_j, Q\in\ca A, q_h\ge 0,~\hbox{for}~h\ne j,~q_j = -1,  |Q|\ge 1\}.
$$
Notice that $Q\in \ca B^+_j$ if and only if we have
$$
\la \widehat \eta^{(k)}, \widehat Q \ra = \eta^{(k)}_j \quad \hbox{for}~k=2, \dots, r\tag{sistemaecalle}
$$
where $\widehat Q = (q_1, \dots, q_{j-1},q_{j+1}, \dots, q_n)\in\N^{n-1}$ and $\widehat \eta^{(k)} = (\eta^{(k)}_1, \dots, \eta^{(k)}_{j-1},\eta^{(k)}_{j+1}, \dots, \eta^{(k)}_n)$, i.e.,~$\widehat Q$ is a solution in $\N^{n-1}$ of the linear system with integer coefficients \rf{sistemaecalle}. 
Moreover, since~$\ca A$ is a submodule of $\Z^n$, Proposition \rf{PrMinComin} applies to $\ca A^+$. Let ${\fr M}= \{M_1,\dots, M_{d}\}$ be the set of minimal elements of $\ca A^+$ and let ${\fr C}=\{C_1,\dots, C_{e}\}$ be the set of cominimal elements of $\ca A^+$ (recall that they all are different from $O$, hence their modulus is at least $1$). We can thus consider the subsets~$\{M'_1,\dots, M'_s\}\subset{\fr M}$ and~$\{C'_1,\dots, C'_{t}\}\subset {\fr C}$ of the minimal and cominimal elements $R$ of~$\ca A^+$ such that $\la \eta^{(1)},  R\ra\in\tau p \Z$.
Then $[\phe]$ can be simplified if and only if there exists $H\in\Z^n$ such that
$$
\la H,M'_h\ra = {1\over\tau p}\la\eta^{(1)},M'_h\ra,
$$
for $1\le h\le s$,
$$
\la H,C'_l\ra = {1\over\tau p}\la\eta^{(1)},C'_{l}\ra,
$$
for $1\le l\le t$, and such that, for any $j=1, \dots, n$, we have
$$
\la\widehat H,\widehat Q\ra - h_j = {1\over\tau p}\left(\la\widehat\eta^{(1)},\widehat Q\ra - \eta^{(1)}_j\right),
$$
for every solution $\widehat Q\in\N^n$ of \rf{sistemaecalle}, with $|Q|\ge 1$, such that $\la \widehat Q, \widehat\eta^{(1)}\ra\in \tau p\, \eta^{(1)}_j\Z$.

%%%%%%%%%%%%%%%%%%%%%%%%%%%%%%%%%%%%%%%%%%%%%%%%%%%%%%%%%%%%%%%%%%%%%%%%%%%%%%%%
\sect Construction of torus actions 
%%%%%%%%%%%%%%%%%%%%%%%%%%%%%%%%%%%%%%%%%%%%%%%%%%%%%%%%%%%%%%%%%%%%%%%%%%%%%%%%

\sm In this last section we shall see some conditions assuring the existence of the torus actions we need.

Let~$X\in\X_n$ be a germ of holomorphic vector field of~$(\C^n,O)$ singular at the origin, in Poincar\'e-Dulac normal form, i.e., 
$$
X= X^{\rm dia} + X^{\rm nil} + X^{\rm res}
$$
where, denoting with $\de_j$ the partial derivative $\de/\de z_j$,
$$
X^{\rm dia} = \sum_{j=1}^n \phe_j z_j \de_j,
$$
$X^{\rm nil}$ is a linear nilpotent vector field singular at the origin such that
$$ 
[X^{\rm dia}, X^{\rm nil}] = 0,
$$
and $X^{\rm res}$ is a holomorphic vector field singular at the origin with no linear part and such that
$$
[X^{\rm dia}, X^{\rm res}] = 0.
$$
In particular
$$ 
[X^{\rm dia}, X^{\rm nil} + X^{\rm res}] = 0.
$$
Recall that the flows of two commuting vector fields also commute (see [Le] Prop. 18.5). We have 
$$
\exp(X^{\rm dia}) = \diag(e^{\phe_1},\dots, e^{\phe_n}) z.
$$
and, in general for a linear vector field $X^{\rm lin} = \sum_{j=1}^n\left(\sum_{h=1}^n a_{hj} z_h\right)\de_j$, we have
$$
\exp(X^{\rm lin}) = e^A z,
$$
where $A$ is the matrix $(a_{hj})$. If $Y$ is a holomorphic vector field singular at the origin with no linear part, then we have
$$
\exp(tY)z= \sum_{k\ge 0} {t^k \over k!}Y^k(z).\tag{esponenziale}
$$
In fact, defining $K_t(z)= z + tY(z)$, we get $K_0(z)=z$ and ${\de\over \de t}K_t(z)|_{t=0} = Y(z)$, then we have~$\exp(tY)z=\lim_{m\to\io} (K_{1/m})^m$, (see [AMR] Theorem $4.1.26$), that is~\rf{esponenziale}.  %[AMR] Theorem $4.1.26$ p.~$254$
Moreover, if~$V, W$ are two commuting vector fields, we have
$$
\exp(t(V+W))=\left( \sum_{k\ge 0} {(tV)^k \over k!}\right)\left( \sum_{k\ge 0}{(tW)^k \over k!}\right)= \exp(tV)\exp(tW).
$$
Then we have the following result. 

\sm\thm{PrFlussi}{Proposition} {\sl Let~$X$ be a germ of holomorphic vector field of~$(\C^n,O)$, singular at the origin, and in Poincar\'e-Dulac normal form. Then its flow is a germ of biholomorphism of~$(\C^n, O)$ in Poincar\'e-Dulac normal form.} 

\sm\proof The flow of~$X^{\rm nil} + X^{\rm res}$ is unipotent, then the linear part of the flow of~$X$ is~$A z$ with~$A$ triangular matrix with diagonal~$\diag(e^{\phe_1},\dots, e^{\phe_n})$, and the flow of $X$ has to commute with the flow of~$X^{\rm dia}$. \qed

\sm In [Zu], Zung found that to find a Poincar\'e-Dulac holomorphic normalization for a germ of holomorphic vector field is the same as to find (and linearize) a suitable torus action which preserves the vector field. To deal with this problem he introduced the notion of {\it toric degree} of a vector field. The following definition is a reformulation of Zung's original one, clearer and more suitable to our needs.

\sm\defin{def1.1} The {\it toric degree} of a germ of holomorphic vector field~$X$ of~$(\C^n, O)$ singular at the origin is the minimum~$r\in\N$ such that the semi-simple part $X^{\rm dia}=\sum_{j=1}^n \phe_j z_j \de_j$ of the linear term of $X$ can be written as linear combination with complex coefficients of $r$ diagonal vector fields with integer coefficients, i.e.,  
$$
X^{\rm dia}= \sum_{k=1}^r \alpha_k Z_k,
$$
where~$\alpha_1, \dots, \alpha_r\in\C^*$ and $Z_k = \sum_{j=1}^n\rho^{(k)}_j z_j \de_j$ with~$\rho^{(k)}\in\Z^n$.
The~$r$-tuple~$Z_1,\dots, Z_k$ is called a {\it $r$-tuple of toric vector fields associated to~$X$}, and the numbers~$\alpha_1, \dots, \alpha_n\in\C$ are a~$r$-tuple of {\it toric coefficients} of the toric $r$-tuple.

\sm In particular, we have
$$
\phe=\sum_{k=1}^r \alpha_k \rho^{(k)}.
$$ 
One can prove  (see [Ra1] pp. 55--57) that the toric coefficients~$\alpha_1, \dots, \alpha_r$ are rationally independent, and $\rho^{(k)}_j = \rho^{(k)}_h$ whenever $\phe_j=\phe_h$, for every $k=1,\dots, r$, implying that 
$$
\res_j^+(\phe)= \bigcap_{k=1}^r\res_j^+(\rho^{(k)})
$$
for any $j=1, \dots, n$.

\sm\thm{ReGrado1CV}{Remark} A vector field has toric degree~$1$ if and only if, chosen a non-zero eigenvalue of its linear part, all the other eigenvalues are rational multiplies of it; then in this case we have uniqueness of the toric vector field associated to $X$ up to multiplication by a non-zero integer.

\sm We recall the following definition from [Zu]

\sm\defin{De4.1} A germ of holomorphic vector field~$X$ of~$(\C^n, O)$ singular at the origin is {\it integrable %in senso non Hamiltoniano
}if it has order $1$ and there exists a positive integer~$1\le m\le n$ such that there exist~$m$ germs of holomorphic vector fields~$X_1=X, X_2, \ldots, X_m$ of~$(\C^n, O)$ singular at the origin and of order $1$, and~$n-m$ germs of holomorphic functions~$g_1, \ldots, g_{n-m}$ in~$(\C^n, O)$ satisfying:{\parindent=30pt
\sm\item{(i)}~$X_1, \ldots, X_m$ commute pairwise and are linearly independent, i.e.,~$X_1\land\cdots\land X_m\not\equiv 0$;
\sm\item{(ii)}~$g_1, \ldots, g_{n-m}$ are common first integrals of~$X_1, \ldots, X_m$, i.e.,~$X_j(g_k)=0$ for any~$j$ and~$k$, and they are functionally independent almost everywhere, i.e.,~$\d g_1\land\cdots\land\d g_{n-m}\not\equiv 0$.
\sm}

Noticing that all the vector fields in the previous definition are integrable, we can define

\sm\defin{De4.3} Let $1\le m\le n$. A {\it set of $m$ integrable vector fields} of $(\C^n, O)$ is a set~$X_1, \ldots, X_m$ of germs of holomorphic vector fields of~$(\C^n, O)$ singular at the origin, of order $1$ and such that:{\parindent=30pt
\sm\item{(i)}~$X_1, \ldots, X_m$ commute pairwise and are linearly independent;
\sm\item{(ii)} there exist~$n-m$ germs of holomorphic functions~$g_1, \ldots, g_{n-m}$ in~$(\C^n, O)$ which are common first integrals of~$X_1, \ldots, X_m$, and they are functionally independent almost everywhere.
\sm}

\sm\thm{Te4.1}{Theorem} (Zung, 2002 [Zu]) {\sl Let~$X$ be a germ of holomorphic vector field of~$(\C^n, O)$ singular at the origin which is integrable. Then~$X$ admits a holomorphic Poincar\'e-Dulac normalization.}

\sm As a corollary of Proposition \rf{PrFlussi}, we obtain

\sm\thm{Co4.1}{Corollary} {\sl The flow of a germ of integrable holomorphic vector field of~$(\C^n, O)$ singular at the origin admits a holomorphic Poincar\'e-Dulac normalization.}

\sm Moreover we have the following result

\sm\thm{Te4.2}{Theorem} (Zung, 2002 [Zu]) {\sl Let $1\le m\le n$. Every set of $m$ integrable vector fields admits a simultaneous holomorphic Poincar\'e-Dulac normalization.}

\sm Thus we have the following corollary

\sm\thm{Co4.1}{Corollary} {\sl Let $1\le m\le n$. The flows of a set of $m$ integrable vector fields admit a simultaneous holomorphic Poincar\'e-Dulac normalization.}

\sm\thm{ReCorPrec}{Remark} The last corollary means that we can conjugate $X_1, \dots, X_m$ to a $m$-tuple of vector fields containing only monomials belonging to the intersection of the additive resonances of the eigenvalues of the linear terms of $X_1, \dots, X_m$. 

\me Now, we introduce an analogous for germs of biholomorphisms of the notion of integrability we described above. 

\sm\defin{De5.1} A germ of biholomorphism~$f$ of~$(\C^n, O)$ fixing the origin {\it commutes with a set of integrable vector fields} if there exists a positive integer~$1\le m\le n$, such that there exists a set of~$m$ germs of holomorphic integrable vector fields~$X_1, \ldots, X_m$ such that 
$$
\d f(X_j)=X_j\circ f
$$
for each~$j=1, \dots, m$.

\sm\thm{ReComm}{Remark} A germ of biholomorphism $f$ of $(\C^n,O)$ commutes with a
vector field $X$ according to the previous definition if and only if it commutes with the flow generated by $X$.

\sm\thm{Te5.2}{Theorem} {\sl Let~$f$ be a germ of biholomorphism of~$(\C^n, O)$ fixing the origin and commuting with a set of integrable holomorphic vector fields $X_1,\dots, X_m$. Then~$f$ commutes with a holomorphic effective action on~$(\C^n,O)$ of a torus of dimension equal to the toric degree $r$ of~$X_1$ and such that the columns of the weight matrix of the action are a~$r$-tuple of toric vectors associated to~$X_1$.}

\sm \proof From the proof of Theorem \rf{Te4.1} (see [Zu]) we get $r$ holomorphic vector field $\ca Z_1, \dots, \ca Z_r$ which generate a $\T^r$ action preserving $X_1,\dots, X_m$. Moreover~$a_{1,k},\dots, a_{m,k}$ holomorphic functions constant on the connected components of each level set $L_y={\bf g}^{-1}({\bf g}(y))$, where we denote by~${\bf g}=(g_1, \dots, g_{n-m})$ the~$(n-m)$-tuple of common first integrals of~$X_1, \dots, X_m$, such that
$$
\ca Z_k = \sum_{j=1}^m a_{j,k} X_j,
$$
for each $k=1,\dots, r$.
Then, for each $k=1,\dots, r$, we have
$$
\eqalign{
\d f({\ca Z}_k) &= \d f\left(\sum_{j=1}^m a_{j,k} X_j\right)\cr
 					 &= \sum_{j=1}^m \left( a_{j,k}\circ f \right)\d f(X_j)\cr
					 &= \sum_{j=1}^m \left (a_{j,k} \circ f\right)\left (X_j\circ f\right)\cr
					 &= \ca Z_k\circ f.}
$$
Thus the torus action commutes with~$f$ as we wanted. \qed

\sm\thm{Co9.1}{Corollary} {\sl Let~$f$ be a germ of biholomorphism of~$(\C^n, O)$ fixing the origin and commuting with a set of integrable holomorphic vector fields. Then~$f$ is holomorphically conjugated to a germ containing only monomials belonging to the intersection of the additive resonances of the eigenvalues of the linear terms of $X_1, \dots, X_m$.}

\sm\proof The assertion follows from Theorem \rf{Te5.2}, Corollary \rf{Co4.1} and Theorem \rf{Te2.1}. \qed

\sm Then we also have the following

\sm\thm{Co9.2}{Corollary} {\sl Let~$f$ be a germ of biholomorphism of~$(\C^n, O)$ fixing the origin and commuting with a set of integrable holomorphic vector fields, such that the intersection of the additive resonances of the eigenvalues of the linear terms of $X_1, \dots, X_m$ is equal or contained in the set of resonances of the spectrum of $\d f_O$. Then~$f$ admits a holomorphic Poincar\'e-Dulac normalization.}

\sm\thm{RemReferee}{Remark} A slight generalization
of the proof of Theorem~\rf{Te4.1} shows that in the statement of Theorem~\rf{Te5.2} it is not necessary for all the vector fields $X_2,\ldots,X_m$ to have order~1; however, it is still necessary that the whole set of vector fields $X_1, \dots, X_m$ commutes with $f$. We refer to [Ra4] for the precise statement and detailed proof.

%%%%%%%%%%%%%%%%%%%%%%%%%%%%%%%%%%%%%%%%%%%%%%%%%%%%%%%%%%%%%%%%%%%%%%%%%%%%%%%%
\vbox{\vskip.75truecm\advance\hsize by 1mm
	\hbox{\centerline{\sezfont References}}
	\vskip.25truecm}\nobreak
%%%%%%%%%%%%%%%%%%%%%%%%%%%%%%%%%%%%%%%%%%%%%%%%%%%%%%%%%%%%%%%%%%%%%%%%%%%%%%%%
\parindent=40pt

\bib{A1} {\sc Abate, M.:} {\sl Discrete local holomorphic dynamics,} in ``Proceedings of 13th Seminar of Analysis and its Applications, Isfahan, 2003'', Eds. S. Azam et al., University of Isfahan, Iran, 2005, pp. 1--32.

\bib{A2} {\sc Abate, M.:} {\sl Discrete holomorphic local dynamical systems,} to appear in ``Holomorphic Dynamical Systems'', G. Gentili, J. Guenot, G. Patrizio eds., Lectures notes in Math., Springer Verlag, 
Berlin, 2009, arXiv:0903.3289v1.

\bib{AMR} {\sc Abraham, R., Marsden, J. E., Ratiu, T.:} ``Manifolds, tensor analysis, and applications'', Second edition, Applied Mathematical Sciences, {\bf 75}, Springer-Verlag, New York, 1988.

\bib{Ar} {\sc Arnold, V. I.:} ``Geometrical methods in the theory of ordinary differen\-tial equations'', Springer-Verlag, Berlin, 1988.

\bib {B} {\sc Bochner, S.:} {\sl Compact Groups of Differentiable Transformation}, Annals of Mathematics (2), {\bf 46}, {\bf 3} (1945), pp. 372--381. 

%\bib{BCS} {\sc Bao, D., Chern, S.-S., Shen, Z.:} ``An introduction to Riemann-Finsler geometry'', Graduate Texts in Mathematics {\bf 200}, Springer-Verlag, New York, 2000.

%\bib{BER} {\sc Baouendi, M. S., Ebenfelt, P., Rothschild, L. P.:} {\sl Dynamics of the Segre varieties of a real submanifold in complex space,} J. Algebraic Geom. {\bf 12} (2003), no. 1, pp. 81--106.

%\bib{BMR} {\sc Baouendi, M. S., Mir, N., Rothschild, L. P.:} {\sl Reflection ideals and mappings between generic submanifolds in complex space,} J. Geom. Anal. {\bf 12} (2002), no. 4, pp. 543--580.

\bib{Bra} {\sc Bracci, F.:} {\sl Local dynamics of holomorphic diffeomorphisms,} Boll. UMI (8), 7--B (2004), pp. 609--636.

\bib{Brj} {\sc Brjuno, A. D.:} {\sl Analytic form of differential equations,} Trans. Moscow Math. Soc., {\bf 25} (1971), pp. 131--288; {\bf 26} (1972), pp. 199--239.

%\bib{C} {\sc Carletti, T.:} {\sl Exponentially long time stability for non-linearizable analytic germs of~$(\C^n, 0)$,} Annales de l'Institut Fourier (Grenoble), {\bf 54} (2004), no. 4, pp. 989--1004.

\bib{D} {\sc Dulac, H.:} {\sl Recherches sur les points singuliers des \'equationes diff\'erentielles}, J. \'Ecole polytechnique II s\'erie cahier IX, (1904), pp. 1--125.

\bib{\'E} {\sc \'Ecalle, J.:} {\sl Singularit\'es non abordables par la g\'eom\'etrie,} Annales de l'Institut Fourier (Grenoble), {\bf 42} (1992), no. 1-2, pp. 73--164.

\bib{\'EV} {\sc \'Ecalle, J., Vallet, B.:} {\sl Correction and linearization of resonant vector fields and diffeomorphisms,} Math. Z., {\bf 229} (1998), no. 2, pp. 249--318.

%\bib{G} {\sc Gray, A.:} {\sl A fixed point theorem for small divisors problems,}  J. Diff. Eq., {\bf 18} (1975), pp. 346--365.

\bib{Le} {\sc J.M. Lee:} ``Introduction to Smooth Manifolds'', Springer-Verlag, New York, 2002.

%\bib{M} {\sc Marmi, S.:} ``An introduction to small divisors problems'', I.E.P.I., Pisa, 2003.

%\bib{N} {\sc Nishimura, Y.:} {\sl Automorphismes analytiques admettant des sous-vari\'et\'es de points fix\'es attractives dans la direction transversale,} J. Math. Kyoto Univ., {\bf 23--2} (1983), pp. 289--299.

\bib{PM} {\sc Perez-Marco, R.}: {\sl Linearization of holomorphic germs with resonant linear part}, Preprint, arXiv:math/0009030v1, 2000.

\bib{Po} {\sc Poincar\'e, H.:} ``\OE uvres, Tome I'', Gauthier-Villars, Paris, 1928, pp. XXXVI--CXXIX.

%\bib{P\"o} {\sc P\"oschel, J.:} {\sl On invariant manifolds of complex analytic mappings near fixed points,} Exp. Math., {\bf 4} (1986), pp. 97--109.

\bib{Ra1} {\sc Raissy, J.}: {\bf Normalizzazione di campi vettoriali olomorfi}. Tesi di Laurea Spe\-cia\-li\-sti\-ca, http://etd.adm.unipi.it/theses/available/etd-06022006-141206/, 2006.

\bib{Ra2} {\sc Raissy, J.}: {\sl Linearization of holomorphic germs with quasi-Brjuno fixed points}, Math. Z., {http://www.springerlink.com/content/3853667627008057/fulltext.pdf}, Online First, (2009).

\bib{Ra3} {\sc Raissy, J.}: {\sl Simultaneous linearization of holomorphic germs in presence of resonances}, to appear in Conform. Geom. Dyn. (2009), arXiv:0812.3579v2.

\bib{Ra4} {\sc Raissy, J.}: {\bf Geometrical methods in the normalization of germs of biholomorphisms}, Ph.D. Thesis, Universit\`a di Pisa, (2009).

%\bib{R\"u} {\sc R\"ussmann, H.:} {\sl Stability of elliptic fixed points of analytic area-preserving mappings under the Brjuno condition,} Ergodic Theory Dynam. Systems, {\bf 22} (2002), pp. 1551--1573.

%\bib{S} {\sc Sternberg, S.:} {\sl Infinite Lie groups and the formal aspects of dynamical systems,} J. Math. Mech., {\bf 10} (1961), pp. 451--474.

\bib{Y} {\sc Yoccoz, J.-C.:} {\sl Th\'eor\`eme de Siegel, nombres de Bruno et polyn\^omes quadratiques,} Ast\'erisque {\bf 231} (1995), pp. 3--88.

\bib{Zu} {\sc Zung, N. T.:} {\sl Convergence versus integrability in Poincar\'e-Dulac normal form}, Math. Res. Lett. {\bf 9}, 2-3, (2002), pp.~217--228.

\bye